\documentclass[10pt, a4paper]{article}
\usepackage{latexsym}
\usepackage{amsmath}
\usepackage{amssymb}
\usepackage[all]{xy}
\usepackage{amsfonts}
\usepackage{verbatim}
\usepackage{amsthm}
\usepackage{mathrsfs}
\usepackage{epsfig}
\usepackage{xy}
\usepackage{array}
\usepackage{stmaryrd}
\usepackage{graphicx,color}
\usepackage{xcolor}
\usepackage{tikz}
\usetikzlibrary{arrows,calc}
\usepackage{etex}
\usepackage{mathdots}
\usepackage{float}
\usepackage{graphics}
\usepackage{pdflscape}

\usepackage{anysize,hyperref}
\input xypic
\xyoption{all}

\usepackage[perpage,symbol]{footmisc}
\usepackage{setspace}
\renewcommand{\cal}{\mathcal}

\def\C{\mathscr{C}}
\def\E{\mathbb{E}}
\def\s{\mathfrak{s}}
\def\id{\mathrm{id}}
\def\op{^\mathrm{op}}
\def\Ab{\mathit{Ab}}
\def\del{\delta}
\def\dr{\ar@{->}[r]}

\def\X{\mathscr{X}}

\def\add{\mbox{add}}
\def\Ext{\mbox{Ext}}
\def\Hom{\mbox{Hom}}

\begin{document}
\baselineskip=15pt
\title{\Large{\bf Triangulated quotient categories revisited\footnote{Supported by the NSF of China (Grants No.\;11131001) and the China Postdoctoral Science Foundation (Grants No.\;2016M591155)}}}
\medskip
\author{\textbf{Panyue Zhou and Bin Zhu}}

\date{}

\maketitle
\def\blue{\color{blue}}
\def\red{\color{red}}

\newtheorem{theorem}{Theorem}[section]
\newtheorem{lemma}[theorem]{Lemma}
\newtheorem{corollary}[theorem]{Corollary}
\newtheorem{proposition}[theorem]{Proposition}
\newtheorem{conjecture}{Conjecture}
\theoremstyle{definition}
\newtheorem{definition}[theorem]{Definition}
\newtheorem{question}[theorem]{Question}
\newtheorem{remark}[theorem]{Remark}
\newtheorem{remark*}[]{Remark}
\newtheorem{example}[theorem]{Example}
\newtheorem{example*}[]{Example}

\newtheorem{construction}[theorem]{Construction}
\newtheorem{construction*}[]{Construction}

\newtheorem{assumption}[theorem]{Assumption}
\newtheorem{assumption*}[]{Assumption}

\baselineskip=17pt
\parindent=0.5cm

\begin{abstract}
\baselineskip=16pt
Extriangulated categories were introduced by Nakaoka and Palu by extracting the similarities between exact categories and triangulated categories. A notion of mutation of subcategories in an extriangulated category is defined in
this article. Let $\cal A$ be an extension closed subcategory of an extriangulated category $\C$. Then the quotient category $\cal M:=\cal A/\X$ carries naturally a triangulated structure whenever $(\cal A,\cal A)$ forms an $\X$-mutation pair. This result unifies many previous constructions of triangulated quotient categories, and using it gives a classification of thick triangulated subcategories of pretriangulated category $\C/\X$, where $\X$ is functorially finite in $\C$. When $\C$ has Auslander-Reiten translation $\tau$, we prove that for a functorially finite subcategory $\X$ of $\C$ containing projectives and injectives, $\C/\X$ is a triangulated category if and only if $(\C,\C)$ is $\X-$mutation, and if and only if  $\tau \underline{\X}=\overline{\X}.$  This generalizes a result by J{\o}rgensen who proved the equivalence between the first and the third conditions for triangulated categories. Furthermore, we show that for such a subcategory $\X$ of the extriangulated category $\C$, $\C$ admits a new extriangulated structure such that $\C$ is a Frobenius extriangulated category. Applications to exact categories and triangulated categories are given. From the applications we present examples that extriangulated categories are neither exact categories nor triangulated categories.\\[0.5cm]
\textbf{Key words:} Extriangulated category; Mutation; Quotient triangulated category\\[0.2cm]
\textbf{ 2010 Mathematics Subject Classification:}  18E10; 18E30; 18E40
\medskip
\end{abstract}

\section{Introduction}

Exact and triangulated categories are two fundamental structures
in algebra and geometry. A first looking at their construction shows that these two kinds categories have many different points:
 while modules or sheaves are forming
abelian exact categories, complexes lead to homotopy or derived categories
that are triangulated. But if we look carefully at the construction inside the categories, they share many similarities: while exact categories admit short exact sequences, triangulated categories admit triangles.
  The similarity between short exact sequences and triangles makes that it is possible to put these two notations into a unified form. This was carried out recently by Nakaoka and Palu [NP].

  After a careful looking what is necessary to define a cotorsion pair on exact categories or triangulated category, Nakaoka and Palu [NP] introduced the notion of an extriangulated category by extracting those properties of $\Ext^1$ on exact categories and on triangulated categories that seem relevant from the point-of-view of cotorsion pairs. Using the frame of extriangulated category, they succeed in unifying and extending the construction on twin cotorsion pairs in triangulated and exact categories. The class of extriangulated categories not only contains exact categories and
extension closed subcategories of triangulated categories as examples, but it is also
closed under taking some quotients. The term ``extriangulated" stands for externally triangulated by means of a bifunctor. It can also be viewed as the mixing of exact and triangulated, or as an
abbreviation of Ext-triangulated.

We want to show further the similarity and difference between exact categories and triangulated categories from the point of view of construction of quotient triangulated categories.

Quotients of exact or triangulated categories have been extensively studied in representation theory. It is in particular very much of interest to know which structure of the quotient inherits from the exact or triangulated category one starts with. Well-studied examples include:
\begin{itemize}
\item stable categories, which are triangulated quotients of Frobenius exact categories. Happel [Ha] shows that if  $(\cal B, \cal S)$ is a Frobenius exact category,
then its stable category $\cal B/\cal I$ carries a triangulated structure, where $\cal I$ is the full subcategory of $\cal B$ consisting of $\cal E$-injective objects. Beligiannis obtained a similar result [B1, Theorem 7.2] by replacing $\cal B$ with a triangulated category $\C$ and replacing $\cal S$ with a proper class of triangles $\cal E$.

\item certain subquotient categories of triangulated categories formed from mutations in triangulated
categories.  Mutation of cluster tilting objects introduced by Buan-Marsh-Reineke-Reiten-Todorov [BMRRT] and by Gei{\ss}-Leclerc-Shr\"{o}r [GLS] is a key integration in attempting to categorify Fomin-Zelevinsky's cluster algebras [FZ] by using quiver representations. The notion of mutation pairs of subcategories in a triangulated category was
defined by Iyama and Yoshino [IY], which is a generalization of mutation of cluster tilting
objects in cluster categories [BMRRT].  We recall the Iyama-Yoshino's definition here.

Let $\cal D\subseteq\cal Z$  be subcategories of a triangulated category $\C$ such that $\Ext^1_{\C}(\cal D,\cal D)=0$. We call $(\cal Z,\cal Z)$ a \emph{$\cal D$-mutation pair} in the sense of Iyama-Yoshino [IY, Definition 2.5], if it satisfies
\begin{itemize}
\item[(1)] For any $X\in\mathcal{Z}$, there exists a triangle $\xymatrix{X\ar[r]^{f}&D\ar[r]^{g}&Y\ar[r]^{h}&X[1]},$ where $Y\in\cal Z$, $f$ is a left $\cal D$-approximation and $g$ is a right $\cal D$-approximation.

\item[(2)] For any $Y\in\mathcal{Z}$, there exists a triangle $\xymatrix{X\ar[r]^{f}&D\ar[r]^{g}&Y\ar[r]^{h}&X[1]},$ where $X\in\cal Z$, $f$ is a left $\cal D$-approximation and $g$ is a right $\cal D$-approximation.
\end{itemize}
Iyama and Yoshino [IY, Theorem 4.2] proved that if $\cal Z$ is extension closed and $(\cal Z,\cal Z)$ is a $\cal D$-mutation pair, then the quotient category $\cal Z/\cal D$ is a triangulated category.  J{\o}rgensen [J, Theorem 2.3] gave a similar construction of triangulated quotient category in another manner. He proved that if $\X$ is a functorially finite subcategory of a triangulated
category $\C$ with Auslander-Reiten translation $\tau$, and if $\X$ satisfies the equation
$\tau \X=\X$, then the quotient category $\C/\X$ is a triangulated category. Later, to unify these two construction of quotient triangulated categories, Liu and Zhu [LZ] defined the notion of $\cal D$-mutation pair of subcategories in a triangulated category without assumption of $\Ext^1_{\C}(\cal D,\cal D)=0$, and showed that Iyama and Yoshino's result.
Recently, Beligiannis [B2] introduced the notion of $\X$-Frobenius subcategory $\cal U$ in a triangulated category $\C$, showed that the quotient category $\cal U/\X$ is a triangulated category. In fact, it is easy to see that Beligiannis' definition is equivalent to Liu-Zhu's definition.

\end{itemize}

The main aim of this article is to give a way to unify the existing different constructions
of triangulated quotient categories and study the properties of the quotients. We define a notion of mutation pairs of subcategories in
an extriangulated category, and show that the corresponding quotient categories carry triangulated structures (see Theorem \ref{h1}). As applications, our result unifies the triangulated quotient category construction considered by  Happel [Ha], Beligiannis [B1, B2], Iyama-Yoshino [IY], J{\o}rgensen [J], and Liu-Zhu [LZ]. We note that there are some other constructions of quotients triangulated categories by using model theory [Li] which is not covered by our construction. An immediate application of our result is to give a classification (see Theorem \ref{d4}) of thick triangulated subcategories of pre-triangulated category $\C/\X$, where $\X$ is a functorially finite in $\C$.
 Let $\X$ be a functorially finite subcategory of a triangulated category $\C$ with the
Auslander-Reiten translation $\tau$.  J{\o}rgensen [J] proved that the quotient category
$\C/\X$ is a triangulated category if and only if $\tau \X=\X$.  We first prove an extriangulated version of J{\o}rgensen theorem and add an equivalent condition that $(\C,\C)$ is an $\X-$mutation (see Theorem \ref{jj1}).
We then prove that for such a subcategory $\X$ of an extriangulated category $\C$, the ambient category $\C$ admits another new extriangulated structure such that $\C$ is a Frobenius extriangulated category with projective subcategory $\X$ (see Theorem \ref{ppp1}). We apply the result to exact categories, we get the exact version of J{\o}rgensen theorem, which completes the results in [KIWY]. If applying to triangulated categories, we not only recover J{\o}rgensen theorem, but also get examples of Frobenius extriangulated categories which are not exact nor triangulated (see Remark \ref{rem1}).

This article is organized as follows: In Section 2, we recall the definition of an extriangulated category and some basic properties which are needed in the proof of our main results. In Section 3, we show certain quotient categories of extriangulated categories form triangulated categories (Theorem \ref{h1}).
In Section 4, we show that an extriangulated version J{\o}rgensen theorem (Theorem \ref{jj1}), and give applications to exact and triangulated categories.

\section{Preliminaries}
Throughout the article, $k$ denotes a field. When we say that $\C$ is a category,
we always assume that $\C$ is a  Krull-Schmidt $k$-linear category.
All the subcategories of a category are full
subcategories and closed under isomorphisms, direct sums and direct summands.
We denote by add$X$ the additive closure generated by object $X$.
We denote by
$\C(A,B)$ or $\Hom_{\C}(A,B)$ the set of morphisms from $A$ to $ B$ in a category $\C$.

We recall some classical terminology which will be needed for our later investigation.
Let $\C$ be a category and $\X$ a subcategory of $\C$. A
morphism $f_B\colon X_B\to B$ of $\C$ with $X_B$ an object in $\X$,
is said to be a \emph{right $\X$-approximation} of $B$, if the
morphism $\C(X,f_B)\colon\C(X,X_B)\to\C(X,B)$ is surjective
for all objects $X$ in $\X$.
Dually, a morphism $g_A\colon A\to X_A$ of $\C$ with $X_A$ an object in $\X$,
is said to be a \emph{left $\X$-approximation} of $C$, if the
morphism $\C(g_A,X_A)\colon\C(X_A,X)\to\C(A,X)$ is surjective
for all objects $X$ in $\X$. Moreover, one says that
$\X$ is \emph{contravariantly finite} in $\C$ if every object in $\C$ has a right $\X$-approximation, \emph{covariantly finite} in $\C$ if every object in $\C$ has a left $\X$-approximation, and \emph{functorially finite} if it is both contravariantly finite and covariantly finite in $\C$. For more details, we refer to [AR].

We recall some basics on extriangulated categories from [NP].

Let $\C$ be an additive category. Suppose that $\C$ is equipped with a biadditive functor $\E\colon\C\op\times\C\to\Ab$. For any pair of objects $A,C\in\C$, an element $\delta\in\E(C,A)$ is called an {\it $\E$-extension}. Thus formally, an $\E$-extension is a triplet $(A,\delta,C)$.
Let $(A,\del,C)$ be an $\E$-extension. Since $\E$ is a bifunctor, for any $a\in\C(A,A')$ and $c\in\C(C',C)$, we have $\E$-extensions
$$ \E(C,a)(\del)\in\E(C,A')\ \ \text{and}\ \ \ \E(c,A)(\del)\in\E(C',A). $$
We abbreviately denote them by $a_\ast\del$ and $c^\ast\del$.
For any $A,C\in\C$, the zero element $0\in\E(C,A)$ is called the spilt $\E$-extension.

\begin{definition}{[NP, Definition 2.3]}\label{a1}
Let $(A,\del,C),(A',\del',C')$ be any pair of $\E$-extensions. A {\it morphism} $$(a,c)\colon(A,\del,C)\to(A',\del',C')$$ of $\E$-extensions is a pair of morphisms $a\in\C(A,A')$ and $c\in\C(C,C')$ in $\C$, satisfying the equality
$$ a_\ast\del=c^\ast\del'. $$
Simply we denote it as $(a,c)\colon\del\to\del'$.
\end{definition}

Let $A,C\in\C$ be any pair of objects. Sequences of morphisms in $\C$
$$\xymatrix@C=0.7cm{A\ar[r]^{x} & B \ar[r]^{y} & C}\ \ \text{and}\ \ \ \xymatrix@C=0.7cm{A\ar[r]^{x'} & B' \ar[r]^{y'} & C}$$
are said to be {\it equivalent} if there exists an isomorphism $b\in\C(B,B')$ which makes the following diagram commutative.
$$\xymatrix{
A \ar[r]^x \ar@{=}[d] & B\ar[r]^y \ar[d]_{\simeq}^{b} & C\ar@{=}[d]&\\
A\ar[r]^{x'} & B' \ar[r]^{y'} & C &}$$

We denote the equivalence class of $\xymatrix@C=0.7cm{A\ar[r]^{x} & B \ar[r]^{y} & C}$ by $[\xymatrix@C=0.7cm{A\ar[r]^{x} & B \ar[r]^{y} & C}]$. For any $A,C\in\C$, we denote as
$ 0=[A\xrightarrow{\binom{1}{0}}A\oplus C\xrightarrow{(0,\ 1)}C].$

\begin{definition}{[NP, Definition 2.9]}\label{a2}
Let $\s$ be a correspondence which associates an equivalence class $\s(\del)=[\xymatrix@C=0.7cm{A\ar[r]^{x} & B \ar[r]^{y} & C}]$ to any $\E$-extension $\del\in\E(C,A)$. This $\s$ is called a {\it realization} of $\E$, if it satisfies the following condition:
\begin{itemize}
\item Let $\del\in\E(C,A)$ and $\del'\in\E(C',A')$ be any pair of $\E$-extensions, with $$\s(\del)=[\xymatrix@C=0.7cm{A\ar[r]^{x} & B \ar[r]^{y} & C}],\ \ \ \s(\del')=[\xymatrix@C=0.7cm{A'\ar[r]^{x'} & B'\ar[r]^{y'} & C'}].$$
Then, for any morphism $(a,c)\colon\del\to\del'$, there exists $b\in\C(B,B')$ which makes the following diagram commutative.
\begin{equation}\label{t1}
\begin{array}{l}
$$\xymatrix{
A \ar[r]^x \ar[d]^a & B\ar[r]^y \ar[d]^{b} & C\ar[d]^c&\\
A'\ar[r]^{x'} & B' \ar[r]^{y'} & C' &}$$
\end{array}
\end{equation}
\end{itemize}
In this case, we say that sequence $\xymatrix@C=0.7cm{A\ar[r]^{x} & B \ar[r]^{y} & C}$ {\it realizes} $\del$, whenever it satisfies $\s(\del)=[\xymatrix@C=0.7cm{A\ar[r]^{x} & B \ar[r]^{y} & C}]$.
Remark that this condition does not depend on the choices of the representatives of the equivalence classes. In the above situation, we say that (\ref{t1}) (or the triplet $(a,b,c)$) {\it realizes} $(a,c)$.
\end{definition}

\begin{definition}{[NP, Definition 2.12]}
We call the pair $(\E,\s)$ an {\it external triangulation} of $\C$ if it satisfies the following conditions:
\begin{itemize}
\item[{\rm (ET1)}] $\E\colon\C\op\times\C\to\Ab$ is a biadditive functor.
\item[{\rm (ET2)}] $\s$ is an additive realization of $\E$.
\item[{\rm (ET3)}] Let $\del\in\E(C,A)$ and $\del'\in\E(C',A')$ be any pair of $\E$-extensions, realized as
$$ \s(\del)=[\xymatrix@C=0.7cm{A\ar[r]^{x} & B \ar[r]^{y} & C}],\ \ \s(\del')=[\xymatrix@C=0.7cm{A'\ar[r]^{x'} & B' \ar[r]^{y'} & C'}]. $$
For any commutative square
$$\xymatrix{
A \ar[r]^x \ar[d]^a & B\ar[r]^y \ar[d]^{b} & C&\\
A'\ar[r]^{x'} & B' \ar[r]^{y'} & C' &}$$
in $\C$, there exists a morphism $(a,c)\colon\del\to\del'$ which is realized by $(a,b,c)$.
\item[{\rm (ET3)$\op$}] Let $\del\in\E(C,A)$ and $\del'\in\E(C',A')$ be any pair of $\E$-extensions, realized by
$$\xymatrix@C=0.7cm{A\ar[r]^{x} & B \ar[r]^{y} & C}\ \ \text{and}\ \ \ \xymatrix@C=0.7cm{A'\ar[r]^{x'} & B' \ar[r]^{y'} & C'}$$
respectively.
For any commutative square
$$\xymatrix{
A \ar[r]^x& B\ar[r]^y \ar[d]^{b} & C\ar[d]^c&\\
A'\ar[r]^{x'} & B' \ar[r]^{y'} & C' &}$$
in $\C$, there exists a morphism $(a,c)\colon\del\to\del'$ which is realized by $(a,b,c)$.

\item[{\rm (ET4)}] Let $(A,\del,D)$ and $(B,\del',F)$ be $\E$-extensions realized by
$$\xymatrix@C=0.7cm{A\ar[r]^{f} & B \ar[r]^{f'} & D}\ \ \text{and}\ \ \ \xymatrix@C=0.7cm{B\ar[r]^{g} & C \ar[r]^{g'} & F}$$
respectively. Then there exist an object $E\in\C$, a commutative diagram
$$\xymatrix{A\ar[r]^{f}\ar@{=}[d]&B\ar[r]^{f'}\ar[d]^{g}&D\ar[d]^{d}\\
A\ar[r]^{h}&C\ar[d]^{g'}\ar[r]^{h'}&E\ar[d]^{e}\\
&F\ar@{=}[r]&F}$$
in $\C$, and an $\E$-extension $\del^{''}\in\E(E,A)$ realized by $\xymatrix@C=0.7cm{A\ar[r]^{h} & C \ar[r]^{h'} & E},$ which satisfy the following compatibilities.
\begin{itemize}
\item[{\rm (i)}] $\xymatrix@C=0.7cm{D\ar[r]^{d} & E \ar[r]^{e} & F}$  realizes $f'_{\ast}\del'$,
\item[{\rm (ii)}] $d^\ast\del''=\del$,

\item[{\rm (iii)}] $f_{\ast}\del''=e^{\ast}\del'$.
\end{itemize}

\item[{\rm (ET4)$\op$}]  Let $(D,\del,B)$ and $(F,\del',C)$ be $\E$-extensions realized by
$$\xymatrix@C=0.7cm{D\ar[r]^{f'} & A \ar[r]^{f} & B}\ \ \text{and}\ \ \ \xymatrix@C=0.7cm{F\ar[r]^{g'} & B \ar[r]^{g} & C}$$
respectively. Then there exist an object $E\in\C$, a commutative diagram
$$\xymatrix{D\ar[r]^{f'}\ar@{=}[d]&E\ar[r]^{f}\ar[d]^{h'}&F\ar[d]^{g'}\\
D\ar[r]^{f'}&A\ar[d]^{h}\ar[r]^{f}&B\ar[d]^{g}\\
&C\ar@{=}[r]&C}$$
in $\C$, and an $\E$-extension $\del^{''}\in\E(C,E)$ realized by $\xymatrix@C=0.7cm{E\ar[r]^{h'} & A \ar[r]^{h} & C},$ which satisfy the following compatibilities.
\begin{itemize}
\item[{\rm (i)}] $\xymatrix@C=0.7cm{D\ar[r]^{d} & E \ar[r]^{e} & F}$  realizes $g'^{\ast}\del$,
\item[{\rm (ii)}] $\del'=e_\ast\del''$,

\item[{\rm (iii)}] $d_\ast\del=g^{\ast}\del''$.
\end{itemize}
\end{itemize}
In this case, we call $\s$ an $\E$-{\it triangulation of }$\C$, and call the triplet $(\C,\E,\s)$ an {\it externally triangulated category}, or for short, {\it extriangulated category} $\C$.
\end{definition}

For an extriangulated category $\C$, we use the following notation:

\begin{itemize}
\item A sequence $\xymatrix@C=0.41cm{A\ar[r]^{x} & B \ar[r]^{y} & C}$ is called a {\it conflation} if it realizes some $\E$-extension $\del\in\E(C,A)$.
\item A morphism $f\in\C(A,B)$ is called an {\it inflation} if it admits some conflation $\xymatrix@C=0.7cm{A\ar[r]^{f} & B \ar[r]& C}.$
\item A morphism $f\in\C(A,B)$ is called a {\it deflation} if it admits some conflation $\xymatrix@C=0.7cm{K\ar[r]& A \ar[r]^f& B}.$
\item If a conflation $\xymatrix@C=0.6cm{A\ar[r]^{x} & B \ar[r]^{y} & C}$ realizes $\del\in\E(C,A)$, we call the pair $(\xymatrix@C=0.41cm{A\ar[r]^{x} & B \ar[r]^{y} & C},\del)$ an {\it $\E$-triangle}, and write it in the following way.
$$\xymatrix{A\ar[r]^{x} & B \ar[r]^{y} & C\ar@{-->}[r]^{\del}&}$$
\item Let $\xymatrix{A\ar[r]^{x}&B\ar[r]^{y}&C\ar@{-->}[r]^{\delta}&}$ and $\xymatrix{A'\ar[r]^{x'}&B'\ar[r]^{y'}&C'\ar@{-->}[r]^{\delta'}&}$ be any pair of $\E$-triangles. If a triplet $(a,b,c)$ realizes $(a,c)\colon\del\to\del'$ as in $(\ref{t1})$, then we write it as
    $$\xymatrix{
A \ar[r]^x \ar[d]_{a} & B\ar[r]^y \ar[d]^{b} & C\ar@{-->}[r]^{\delta} \ar[d]^{c}&\\
A'\ar[r]^{x'} & B' \ar[r]^{y'} & C' \ar@{-->}[r]^{\delta'}&}$$
and call $(a,b,c)$ a {\it morphism of $\E$-triangles}.
\end{itemize}

\begin{remark}
For any object $A,B\in\C$, we have  $\E$-triangles in $\C$: $$\xymatrix{A\ar[r]^{\binom{1}{0}\quad}&A\oplus B\ar[r]^{\quad(0,\ 1)}&B\ar@{-->}[r]^{0}&}\ \ and \ \ \xymatrix{A\ar[r]^{\binom{0}{1}\quad}&B\oplus A\ar[r]^{\quad(1,\, 0)}&B\ar@{-->}[r]^{0}&}.$$
\end{remark}

\begin{example}
(1) Exact category $\cal B$ can be viewed as an extriangulated category. For the definition and basic properties of an exact category, see [B\"{u}] and [Ke]. In fact, a biadditive functor $\E:=\Ext^1_{\cal B}\colon
\cal B\op\times\cal B\to \Ab$.  Let $A,C\in\cal B$ be any pair of objects. Define $\Ext^1_{\cal B}(C,A)$ to be the collection of all equivalence classes of short exact sequences of the form $\xymatrix@C=0.7cm{A\ar[r]^{x} & B \ar[r]^{y} & C}$. We denote the equivalence class by $[\xymatrix@C=0.7cm{A\ar[r]^{x} & B \ar[r]^{y} & C}]$ as before.
For any $\del=[\xymatrix@C=0.7cm{A\ar[r]^{x} & B \ar[r]^{y} & C}]\in\Ext^1_{\cal B}(C,A)$,
define the realization $\s(\del)$ of $[\xymatrix@C=0.7cm{A\ar[r]^{x} & B \ar[r]^{y} & C}]$ to be $\del$ itself.  For more details, see [NP, Example 2.13].
\vspace{2mm}

(2) Let $\C$ be an triangulated category with shift functor $[1]$.
Put $\E:=\C(-,-[1])$. For any $\delta\in\E(C,A)=\C(C,A[1])$, take a triangle
$$\xymatrix@C=0.7cm{A\ar[r]^{x} & B \ar[r]^{y} & C\ar[r]^{\del}&A[1]}$$
and define as $\s(\del)=[\xymatrix@C=0.7cm{A\ar[r]^{x} & B \ar[r]^{y} & C}]$. Then $(\C,\E,\s)$
is an extriangulated category. It is easy to see that extension closed subcategories of triangulated categories are also
extriangulated categories. For more details, see [NP, Proposition 3.22].
\vspace{2mm}

(3) Let $\C$ be an extriangulated category, and $\cal J$ a subcategory of $\C$.
If $\cal J\subseteq\cal P\cap\cal I$, where $\cal P$ is the full category of projective objects in $\C$ and $\cal I$ is the full category of injective objects in $\C$, then $\C/\cal J$ is an extriangulated category.
This construction gives extriangulated categories which are not exact nor triangulated in general.
For more details, see [NP, Proposition 3.30]. We will give many examples of extriangulated categories which are not exact nor triangulated categories.
\end{example}

There are some basic results on extriangulated categories which are needed later on.

\begin{lemma}\label{a6}
Let $\C$ be an extriangulated category, $$\xymatrix{A\ar[r]^{x}&B\ar[r]^{y}&C\ar@{-->}[r]^{\delta}&}$$
an $\E$-triangle. Then we have the following long exact sequence:
$$\C(-, A)\xrightarrow{\C(-,x)}\C(-, B)\xrightarrow{\C(-,y)}\C(-, C)\xrightarrow{\delta^{\sharp}_-}
\E(-, A)\xrightarrow{\E(-,x)}\E(-, B)\xrightarrow{\E(-,y)}\E(-, C)
,$$
$$\C(C,-)\xrightarrow{\C(y,-)}\C(B,-)\xrightarrow{\C(x,-)}\C(A,-)\xrightarrow{\delta_{\sharp}^-}
\E(C,-)\xrightarrow{\E(y,-)}\E(B,-)\xrightarrow{\E(x,-)}\E(A,-)
.$$
\end{lemma}

\proof See the proofs of Proposition 3.3 and Proposition 3.13 in [NP].

\begin{lemma}\label{a3}{\emph{[NP, Corollary 3.5]}}
Let $(\C,\E,\s)$ be a extriangulated category, and
$$\xymatrix{
A \ar[r]^x \ar[d]_{a} & B\ar[r]^y \ar[d]^{b} & C\ar@{-->}[r]^{\delta} \ar[d]^{c}&\\
A'\ar[r]^{x'} & B' \ar[r]^{y'} & C' \ar@{-->}[r]^{\delta'}&}$$
 any morphism of $\E$-triangles. Then the following are equivalent.

\emph{(1)}  $a$ factors through $x$;

\emph{(2)}  $a_{\ast}\delta=c^\ast\delta'=0$;

\emph{(3)}  $c$ factors through $y'$.\\
In particular, in the case $\delta=\delta'$ and $(a,b,c)=(\id,\id,\id)$, we obtain
$$ x\ \text{is a section}\ \ \Leftrightarrow\ \delta\ \text{splits}\ \Leftrightarrow\ y\ \text{is a retraction}. $$
\end{lemma}

\begin{lemma}\label{c6}{\emph{[NP, Corollary 3.6]}}
Given a morphism of $\E$-triangles
$$\xymatrix{A\ar[r]^{x}\ar[d]^a&B\ar[r]^y\ar[d]^b&C\ar@{-->}[r]^\del\ar[d]^c&\\
A'\ar[r]^{x'}&B'\ar[r]^{y'}&C'\ar@{-->}[r]^{\del'}&}
$$
in $\C$.\\
\emph{(1)} If $a,b$ are isomorphisms, then $c$ is an isomorphism;\\
\emph{(2)} If $a,c$ are isomorphisms, then $b$ is an isomorphism;\\
\emph{(3)} If $b,c$ are isomorphisms, then $a$ is an isomorphism.
\end{lemma}

\begin{lemma}\label{c3}
 Let $x\colon A\to B$, $y\colon D\to C$ and $f\colon A\to C$ be any morphisms in $\C$. \\
\emph{ (1) }If $x$ is an inflation, then
$\binom{f}{x}\colon A\to C\oplus B$ is an inflation in $\C$.\\
\emph{ (2)} If $y$ is an deflation, then
$(y,f)\colon D\oplus A\to C$ is an inflation in $\C$.
\end{lemma}

\proof  See the proof of Corollary 3.16 in [NP].

\section{Triangulated quotient categories}
\setcounter{equation}{0}

Given two objects $A$ and $B$ of a category $\C$, we denote by $\X(A,B)$ the set of morphisms from $A$ to $B$ of $\C$ which
factor through an object of $\X$. It is well known that $\X(A,B)$ is a subgroup of $\C(A,B)$, and that the family of
these subgroups $\X(A,B)$ forms an ideal of $\C$. Thus we
have the category $\C/\X$ whose objects are objects of
$\C$ and whose morphisms are elements of $\C(A,B)/\X(A,B)$. The composition of $\C/\X$ is induced
canonically by the composition of $\C$. We then have an additive functor $\pi\colon \C\to \C/\X$,
$\pi(A)=\overline{A}$ and $\pi(f)=\overline{f}$. Such category $\C/\X$ is called the \emph{quotient category} of $\C$
by $\X$.

\begin{definition}\label{b1}
Let $\X$ be a subcategory of an extriangulated category $\C$. A morphism $f\colon A\to B$ in $\C$ is called
\emph{$\X$-monic}, if for any object $X\in\X$, we have that
$$\C(B,X)\xrightarrow{\C(f,X)}\C(A,X)\to 0$$ is exact.
Dually, a morphism $g\colon C\to D$ in $\C$ is called
\emph{$\X$-epic}, if for any object $X\in\X$, we have that $$
\C(X,C)\xrightarrow{\C(X,g)}\C(X,D)\to 0$$ Obviously,
any left $\X$-approximation is $\X$-monic and any right
$\X$-approximation is $\X$-epic.
\end{definition}

\begin{definition}\label{b2}
Let $\X, \cal A$ and $\cal B$ be subcategories of an extriangulated category $\C$, and $\X\subseteq\cal A$ and $\X\subseteq\cal B$. The pair $(\cal A,\cal B)$ is called an $\X$-\emph{mutation pair }if it satisfies
\begin{itemize}
\item[$(1)$] For any $A\in\cal A$, there exists an  $\E$-triangle
$$\xymatrix{A\ar[r]^{\alpha}&X\ar[r]^{\beta}&B\ar@{-->}[r]^{\delta}&}$$
where $B\in\cal B$, $\alpha$ is a left $\X$-approximation of $A$ and $\beta$ is a right $\X$-approximation of $B$.

\item[$(2)$] For any $B\in\cal B$, there exists an $\E$-triangle
$$\xymatrix{A\ar[r]^{\alpha}&X\ar[r]^{\beta}&B\ar@{-->}[r]^{\delta}&}$$
where $A\in\cal A$, $\alpha$ is a left $\X$-approximation of $A$ and $\beta$ is a right $\X$-approximation of $B$.
\end{itemize}
Note that if $\C$ is a triangulated category, $\X$-mutation pair is just the same as Liu-Zhu's definition
[LZ, Definition 2.6].
In addition, if $\X$ is a rigid subcategory of $\C$, i.e., $\Ext_{\C}^1(\X,\X)=0$, then $\X$-mutation pair is just the same as Iyama-Yoshino's definition [IY, Definition 2.5].
\end{definition}

\begin{definition}\label{b3}
Let $\C$ be an extriangulated category. A subcategory $\cal A$ of $\C$
is called \emph{extension closed }if for any an $\E$-triangle in $\C$:
$\xymatrix{A\ar[r]^{x}&B\ar[r]^{y}&C\ar@{-->}[r]^{\delta}&,}$
 $A,C\in\cal A$ implies $B\in\cal A$.
\end{definition}

\begin{definition}\label{e3}
Let $\C$ be an extriangulated category $\C$. A subcategory $\X$ of $\C$ is called
\emph{rigid} subcategory of $\C$, if $\E(\X,\X)=0$, i.e. $\E(A,B)=0$, for any $A,B\in\X$.
\end{definition}

\begin{lemma}\label{e4}
Assume that  $(\cal A,\cal B)$ forms an $\X$-mutation pair. If $\X$ is a rigid subcategory of an extriangulated category $\C$. Then
$$\E(\X,\cal A)=0\;\,and\;\, \E(\cal B,\X)=0.$$
\end{lemma}

\proof For any $A\in\cal A$, there exists an $\E$-triangle
$$\xymatrix{A\ar[r]^{f}&X_A\ar[r]^{g}&B\ar@{-->}[r]^{\del}&,}$$
where $X_A\in\X, B\in\cal B$, $f$ is a left $\X$-approximation and $g$ is a right $\X$-approximation.
Applying the functor $\C(\X,-)$ to this $\E$-triangle, we have the following exact sequence
$$\C(\X,X_A)\xrightarrow{\C(\X,g)}\C(\X,B)\xrightarrow{}\E(\X,A)\xrightarrow{}\E(\X,X_A)=0.$$
It follows that $\E(\X,A)=0$ for any $A\in\cal A$ and then $\E(\X,\cal A)=0$.
\vspace{2mm}

Dually, we can prove $\E(\cal B,\X)=0$.  \qed

\subsection{Basics on quotient categories}
In this subsection, let $\C$ be an extriangulated category and $\X\subseteq\cal A$ subcategories of $\C$.

\begin{lemma}\label{b4}
Let
$$\xymatrix{A\ar[r]^x\ar[d]^{a}&B\ar[r]^y\ar[d]^{b}&C\ar@{-->}[r]^{\delta}\ar[d]^{c}&\\
A'\ar[r]^{x'}&X'\ar[r]^{y'}&C'\ar@{-->}[r]^{\delta'}&}$$
be morphisms of $\E$-triangles in $\C$, where $X'\in\X$ and $x$ is $\X$-monic. Then $\overline{a}=0$ in quotient category $\cal A/\X$ implies $\overline{c}=0$
in quotient category $\cal A/\X$.
\end{lemma}

\proof
Since $\overline{a}=0$, there exist $X\in\X$, $a_1\colon A\rightarrow X$ and $a_2\colon X\rightarrow B$ such that
$a=a_2a_1$. Since $X\in\X$ and $x$ is $\X$-monic, there exists $a_3\colon B\rightarrow X$ such that $a_3x=a_1$.
It follows that $a=a_2a_1=a_2a_3x$, that is to say, $a$ factors through $x$.
 By Lemma \ref{a3}, we have that $c$ factors through $y'$. It follows that
$\overline{c}=0$. \qed

\begin{lemma}\label{b5}
Let $\xymatrix{A\ar[r]^{x}&X\ar[r]^{y}&C\ar@{-->}[r]^{\delta}&}$ and $\xymatrix{A\ar[r]^{x'}&X'\ar[r]^{y'}&C'\ar@{-->}[r]^{\delta'}&}$
be $\E$-triangles in $\C$, where $x,x'$ are left $\X$-approximations. Then $C$ and $C'$ are isomorphic in $\cal A/\X$.
\end{lemma}

\proof Since $x,x'$ are left $\X$-approximations, we obtain a morphism of $\E$-triangles
$$\xymatrix{A\ar[r]^x\ar@{=}[d]&X\ar[r]^y\ar[d]^{b}&C\ar@{-->}[r]^{\delta}\ar[d]^{c}&\\
A\ar[r]^{x'}\ar@{=}[d]&X'\ar[r]^{y'}\ar[d]^{b'}&C\ar@{-->}[r]^{\delta'}\ar[d]^{c'}&\\
A\ar[r]^{x}&X\ar[r]^{y}&C\ar@{-->}[r]^{\delta}&}$$
By Lemma \ref{b4}, we have that $\overline{c'}\circ\overline{c}=\overline{\id_C}$.
Similarly, we can show that $\overline{c}\circ\overline{c'}=\overline{\id_C'}$.
Therefore, $C\simeq C'$ in $\cal A/\X$. \qed

\begin{construction}\label{b6}
Assume that $(\cal A,\cal B)$ is an $\X$-mutation pair. For any object $A\in\cal A$, take an $\E$-triangle
$\xymatrix{A\ar[r]^{\alpha}&X\ar[r]^{\beta}&B\ar@{-->}[r]^{\delta}&,}$
where $X\in\X, B\in\cal B$, $\alpha$ is a left $\X$-approximation and $\beta$ is a right $\X$-approximation.

Define $\mathbb{G}(A)=\mathbb{G}A$ to be the image of $B$ in $\cal B/\X$.
This is well-defined by Lemma \ref{b5}.

For any morphism $a\in\cal A(X,X')$, since $\alpha$ is a left $\X$-approximation and $X'\in\X$, there exists
a morphism $b\colon X\to X'$, such that $b\alpha=\alpha'a$. By (ET3), there exists a morphism $c\colon B\to B'$ such that $(a,b,c)$ is a morphism of $\E$-triangles.
$$\xymatrix{A\ar[r]^{\alpha}\ar[d]^a&X\ar[r]^{\beta}\ar[d]^{b}&B\ar@{-->}[r]^\del\ar[d]^{c}&\\
A'\ar[r]^{\alpha'}&X'\ar[r]^{\beta'}&B'\ar@{-->}[r]^{\del'}&}$$
For any $\overline{a}\in\cal A/\X(A,A')$, define $\mathbb{G}\overline{a}$ to be the image
$\overline{c}$ of $c$ in $\cal B/\X$. This is well-defined by Lemma \ref{b4}, and the following proposition holds.
\end{construction}

\begin{proposition}\label{b7}
The functor $\mathbb{G}\colon\cal A/\X\rightarrow\cal B/\X$ is an equivalence.
\end{proposition}

\proof For any $B\in\cal B$, we take an $\E$-triangle
$$\xymatrix{\mathbb{H}B\ar[r]^{\mu}&X_0\ar[r]^{\nu}&B\ar@{-->}[r]^{\eta}&,}$$
where $X_0\in\X, \mathbb{H}B\in\cal A$, $\mu$ is a left $\X$-approximation and $\nu$ is a right $\X$-approximation. We can construct a functor $\mathbb{H}\colon\cal B/\X\rightarrow\cal A/\X$ in a dual manner. We can easily show that $\mathbb{H}$ gives a quasi-inverse of $\mathbb{G}$.  \qed

\subsection{Mutaions and triangulated quotient categories}
\setcounter{equation}{0}
In this subsection, let $\C$ be an extriangulated category and let $\X\subseteq\cal A$ be subcategories
of $\C$. We assume that the following two conditions concerning $\X$ and $\cal A$:

(1) $\cal A$ is extension closed;

(2) $(\cal A,\cal A)$ forms an $\X$-mutation pair.\\
Under such a setting, we put
$$\cal M:=\cal A /\X.$$
The aim of this section is to prove that $\cal M$ admits a triangulated structure.

\begin{definition}\label{c1}
Let
$$\langle1\rangle:=\mathbb{G}\colon\cal M\rightarrow\cal M$$
be the equivalence constructed in Proposition \ref{b7}. Thus for any $A\in\cal A$, we have an $\E$-triangle
$$\xymatrix{A\ar[r]^{\alpha}&X\ar[r]^{\beta}&A\langle1\rangle\ar@{-->}[r]^{\delta}&,}$$
where $X\in\X, A\langle1\rangle\in\cal A$, $\alpha$ is a left $\X$-approximation and $\beta$ is a right $\X$-approximation.

Let
$$\xymatrix{A\ar[r]^{x}&B\ar[r]^{y}&C\ar@{-->}[r]^\del&}$$
 be an $\E$-triangle in $\C$ with $A,B,C\in\cal A$ and $x$ is $\X$-monic. Since $x$ is $\X$-monic,
there exists the following commutative diagram:
\begin{equation}\label{t2}
\begin{array}{l}
\xymatrix{A\ar[r]^{x}\ar@{=}[d]&B\ar[r]^{y}\ar[d]^{b}&C\ar@{-->}[r]^\del\ar[d]^{z}&\\
A\ar[r]^{\alpha}&X\ar[r]^{\beta}&A\langle1\rangle\ar@{-->}[r]^{\del'}&.}
\end{array}
\end{equation}
Then we have the following sextuple in the quotient category $\cal M$:
\begin{equation}\label{t3}
\xymatrix{A\ar[r]^{\overline{x}}&B\ar[r]^{\overline{y}}&C\ar[r]^{\overline{z}\quad}&A\langle1\rangle}.
\end{equation}
We define the \emph{standard triangles} in $\cal M$ as the sextuples which are isomorphic to (\ref{t3}).
\end{definition}

\begin{lemma}\label{c2}
Given a morphism of $\E$-triangles in $\C$
\begin{equation}\label{t5}
\begin{array}{l}
\xymatrix{A\ar[r]^{x}\ar[d]^a&B\ar[r]^y\ar[d]^b&C\ar@{-->}[r]^\del\ar[d]^c&\\
A'\ar[r]^{x'}&B'\ar[r]^{y'}&C'\ar@{-->}[r]^{\del'}&,}
\end{array}
\end{equation}
where $x$ and $x'$ are $\X$-monic and $A,B,C,A',B',C'\in\cal A$. Then we have the morphism of standard triangles in $\cal M$:
$$\xymatrix{A\ar[r]^{\overline{x}}\ar[d]^{\overline{a}}&B\ar[r]^{\overline{y}}\ar[d]^{\overline{b}}&C\ar[r]^{\overline{z}\quad}\ar[d]^{\overline{c}}& A\langle1\rangle\ar[d]^{\overline{a}\langle1\rangle}\\
A'\ar[r]^{\overline{x'}}&B'\ar[r]^{\overline{y'}}&C'\ar[r]^{\overline{z'}\quad}&
A'\langle1\rangle}$$
\end{lemma}

\proof Consider the following commutative diagrams where the rows are $\E$-triangles in $\C$, and $\overline{a'}=\overline{a}\langle1\rangle$:
\begin{equation}\label{t6}
\begin{array}{l}
\xymatrix{A\ar[r]^{x}\ar@{=}[d]&B\ar[r]^y\ar[d]^d&C\ar@{-->}[r]^\del\ar[d]^z&\\
A\ar[r]^{\alpha_A}&X\ar[r]^{\beta_A}&A\langle1\rangle\ar@{-->}[r]^{\eta}&,}
\end{array}
\end{equation}

\begin{equation}\label{t7}
\begin{array}{l}
\xymatrix{A\ar[r]^{\alpha_A}\ar[d]^a&X\ar[r]^{\beta_A}\ar[d]^f&A\langle1\rangle\ar@{-->}[r]^\eta\ar[d]^{a'}&\\
A'\ar[r]^{\alpha_{A'}}&X'\ar[r]^{\beta_{A'}}&A'\langle1\rangle\ar@{-->}[r]^{\eta'}&,}
\end{array}
\end{equation}
and
\begin{equation}\label{t8}
\begin{array}{l}
\xymatrix{A'\ar[r]^{x'}\ar@{=}[d]&B'\ar[r]^{y'}\ar[d]^g&C'\ar@{-->}[r]^{\del'}\ar[d]^{z'}&\\
A'\ar[r]^{\alpha_A}&X'\ar[r]^{\beta_A}&A'\langle1\rangle\ar@{-->}[r]^{\eta'}&.}
\end{array}
\end{equation}
By composing the commutative diagram (\ref{t6}) and (\ref{t7}), we have the following commutative diagram:
$$\xymatrix{A\ar[r]^{x}\ar[d]^{a}&B\ar[r]^y\ar[d]^{fd}&C\ar@{-->}[r]^\del\ar[d]^{a'z}&\\
A'\ar[r]^{\alpha_{A'}}&X'\ar[r]^{\beta_{A'}}&A'\langle1\rangle\ar@{-->}[r]^{\eta'}&.}$$
By composing the commutative diagram (\ref{t5}) and (\ref{t8}), we have the following commutative diagram:
$$\xymatrix{A\ar[r]^{x}\ar[d]^a&B\ar[r]^y\ar[d]^{gb}&C\ar@{-->}[r]^\del\ar[d]^{z'c}&\\
A'\ar[r]^{\alpha_{A'}}&X'\ar[r]^{\beta_{A'}}&A'\langle1\rangle\ar@{-->}[r]^{\eta'}&.}$$
By Lemma \ref{b4}, we have $\overline{a}\langle1\rangle\cdot\overline{z}=\overline{z'}\cdot\overline{c}$. \qed

\begin{lemma}\label{c4}
Assume that we are given $\E$-triangles in $\C$:
$\xymatrix{A\ar[r]^{f}&B\ar[r]^{f'}&C\ar@{-->}[r]^{\del_f}&},$
$$\xymatrix{A\ar[r]^{g}&D\ar[r]^{g'}&F\ar@{-->}[r]^{\del_g}&} \  and  \
\xymatrix{E\ar[r]^{h}&B\ar[r]^{h'}&D\ar@{-->}[r]^{\del_h}&}$$
satisfying $h'f=g$. Then there exists an $\E$-triangle
$\xymatrix{E\ar[r]^{x}&C\ar[r]^{y}&F\ar@{-->}[r]^{\del}&}$
which makes
$$\xymatrix{&E\ar@{=}[r]\ar[d]^{h}&E\ar[d]^{x}\\
A\ar[r]^{f}\ar@{=}[d]&B\ar[r]^{f'}\ar[d]^{h'}&C\ar@{-->}[r]^{\delta_f}\ar[d]^{y}&\\
A\ar[r]^{g}&D\ar[r]^{g'}\ar@{-->}[d]^{\delta_h}&F\ar@{-->}[r]^{\delta_g}\ar@{-->}[d]^{\del}&\\
&&}$$
commutative in $\C$.
\end{lemma}

\proof See the proof the dual of Proposition 3.17 in [NP].

\begin{remark}\label{c5}
In Lemma \ref{c4}, if $g$ and $h$ are $\X$-monic, we have that $x$ is $\X$-monic.
\end{remark}

\proof Let $u\colon E\to X$ be any morphism in $\C$, where $X\in\X$. Since $h$ is $\X$-monic, there exists a morphism $v\colon B\to X$ such that $u=vh$. Since $g$ is $\X$-monic, there exists a morphism $w\colon D\to X$ such that $wg=vf$. Note that $(v-wh')f=vf-wg=0$. By Lemma \ref{a6}, there exists a morphism $s\colon C\to X$ such that $v-wh'=sf'$. It follows that $u=vh=sf'h=sx$. \qed

\begin{lemma}\label{c7}
Let $\xymatrix{A\ar[r]^{x}&B\ar[r]^{y}&C\ar@{-->}[r]^{\del}&}$ and $\xymatrix{A\ar[r]^{\binom{x}{\alpha}\quad}&B\oplus X\ar[r]^{\quad(m,\ \gamma)}&B\ar@{-->}[r]^{\eta}&}$ be two $\E$-triangles in $\cal A$. If $x$ is $\X$-monic and $\alpha$ is a left $\X$-approximation of $A$, then
we have an isomorphism of standard triangles in $\cal M$:
$$\xymatrix{A\ar[r]^{\overline{x}}\ar@{=}[d]&B\ar[r]^{\overline{m}}\ar@{=}[d]&D\ar[r]^{\overline{n}\quad}\ar[d]^{\overline{v}}& A\langle1\rangle\ar@{=}[d]\\
A\ar[r]^{\overline{x}}&B\ar[r]^{\overline{y}}&C\ar[r]^{\overline{z}\quad}&
A\langle1\rangle.}$$
\end{lemma}

\proof By Lemma \ref{c4}, we obtain a commutative diagram
$$\xymatrix{&X\ar@{=}[r]\ar[d]^{\binom{0}{1}}&X\ar[d]^{\gamma}\\
A\ar[r]^{\binom{x}{\alpha}\quad}\ar@{=}[d]&B\oplus X\ar[r]^{\quad(m,\ \gamma)}\ar[d]^{(1,\ 0)}&D\ar@{-->}[r]^{\eta}\ar[d]^{v}&\\
A\ar[r]^{x}&B\ar[r]^{y}\ar@{-->}[d]^{0}&C\ar@{-->}[r]^{\delta}\ar@{-->}[d]^{\varepsilon}&\\
&&}$$
of $\E$-triangles. By Lemma \ref{c2}, we have the morphism of standard triangles in $\cal M$:
$$\xymatrix{A\ar[r]^{\overline{x}}\ar@{=}[d]&B\ar[r]^{\overline{m}}\ar@{=}[d]&D\ar[r]^{\overline{n}\quad}\ar[d]^{\overline{v}}& A\langle1\rangle\ar@{=}[d]\\
A\ar[r]^{\overline{x}}&B\ar[r]^{\overline{y}}&C\ar[r]^{\overline{z}\quad}&
A\langle1\rangle.}$$
It remains to show that $\overline{v}$ is an isomorphism in $\cal M$.
Since $x$ is $\X$-monic and $X\in\X$, there exists a morphism $w\colon B\to X$ such that $\alpha=wx$.
Thus, we obtain a commutative diagram
$$\xymatrix{A\ar[r]^{x}\ar@{=}[d]&B\ar[r]^{y}\ar[d]^{\binom{1}{w}}&C\ar@{-->}[r]^{\delta}\ar[d]^{v'}&\\
A\ar[r]^{\binom{x}{\alpha}\quad}\ar@{=}[d]&B\oplus X\ar[r]^{\quad(m,\gamma)}\ar[d]^{(1,\, 0)}&D\ar@{-->}[r]^{\eta}\ar[d]^{v}&\\
A\ar[r]^{x}&B\ar[r]^{y}&C\ar@{-->}[r]^{\delta}&}$$
of $\E$-triangles.
By Lemma \ref{c6}, we have that $vv'$ is an isomorphism and then $\overline{v}\circ\overline{v'}$
is an isomorphism in $\cal M$.
On the other hand, since $(w,\,1)\binom{x}{\alpha}=0$, by Lemma \ref{a6}, there exists a morphism
$s\colon D\to X$ such that $s(m,\, \gamma)=(w,\, -1)$. In particular, $sm=w$ and $s\gamma=-1$.
Note that $(1+\gamma s)\gamma=\gamma+\gamma s\gamma=\gamma-\gamma=0$. By Lemma \ref{a6}, there exists a morphism $t\colon C\to D$ such that $1+\gamma s=tv$. It follows that $\overline{t}\circ\overline{v}=\overline{1}$. Therefore, $\overline{v}$ is an isomorphism in $\cal M$. \qed

Now we can state and show our main theorem in this section.

\begin{theorem}\label{h1}
 The quotient category $\cal M:=\cal A/\X$ is a triangulated category with respect to the auto-equivalence $\langle1\rangle$ and standard triangles defined in Definition \emph{\ref{c1}}.
\end{theorem}

\proof We will check the axioms of triangulated categories.
\vspace{1mm}

(TR1) By definition, the class of standard triangles is closed under isomorphisms. The following commutative
diagram
$$\xymatrix{A\ar[r]^{\id_A}\ar@{=}[d]&A\ar[r]^{0}\ar[d]^{\alpha}&0\ar@{-->}[r]^\del\ar[d]^{0}&\\
A\ar[r]^{\alpha}&X\ar[r]^{\beta}&A\langle1\rangle\ar@{-->}[r]^{\del'}&.}$$
shows that $\xymatrix{A\ar[r]^{\overline{\id_A}}&A\ar[r]^{\overline{0}}&0\ar[r]^{\overline{0}\quad}&A\langle1\rangle}$
is a standard triangle.

Let $x\in\cal M(A,B)$ be any morphism. For object $A\in\cal A$, there exists an $\E$-triangle in $\C$:
$$\xymatrix{A\ar[r]^{\alpha_A}&X\ar[r]^{\beta_A\;}&A\langle1\rangle\ar@{-->}[r]^{\eta}&.}$$
By Lemma \ref{c3}, we obtain an $\E$-triangle in $\C$:
$\xymatrix{A\ar[r]^{\binom{x}{\alpha}\quad\;}&B\oplus X\ar[r]^{\quad(y,\ \gamma)}&C\ar@{-->}[r]^{\delta}&.}$
By Lemma \ref{c4}, we have a commutative diagram
$$\xymatrix{&B\ar@{=}[r]\ar[d]^{\binom{1}{0}}&B\ar[d]^{u}\\
A\ar[r]^{\binom{x}{\alpha_A}\quad}\ar@{=}[d]&B\oplus X\ar[r]^{\quad(y,\ \gamma)}\ar[d]^{(0,\ 1)}&C\ar@{-->}[r]^{\delta}\ar[d]^{v}&\\
A\ar[r]^{\alpha_A}&X\ar[r]^{\beta_A}\ar@{-->}[d]^{0}&A\langle1\rangle\ar@{-->}[r]^{\eta}\ar@{-->}[d]^{\varepsilon}&\\
&&}$$
of $\E$-triangles in $\C$. Since $\cal A$ is extension closed, we have $C\in\cal A$.
Thus we have a standard triangle
$\xymatrix{A\ar[r]^{\overline{x}}&B\ar[r]^{\overline{y}}&C\ar[r]^{\overline{v}\quad}&A\langle1\rangle}$.
\vspace{1mm}

(TR2) Let $\xymatrix{A\ar[r]^{\overline{x}}&B\ar[r]^{\overline{y}}&C\ar[r]^{\overline{z}\quad}&A\langle1\rangle}$
be a standard triangle in $\cal M$. By Lemma \ref{c7}, we may assume that it is induced by the $\E$-triangle
$$\xymatrix{A\ar[r]^{\binom{x}{\alpha_A}\quad\;}&B\oplus X_A\ar[r]^{\quad(y,\ \gamma)}&C\ar@{-->}[r]^{\delta}&}.$$
By Lemma \ref{c4}, we obtain a commutative diagram
\begin{equation}\label{q1}
\begin{array}{l}
\xymatrix{&B\ar@{=}[r]\ar[d]^{\binom{1}{0}}&B\ar[d]^{y}\\
A\ar[r]^{\binom{x}{\alpha_A}\quad}\ar@{=}[d]&B\oplus X_A\ar[r]^{\quad(y,\ \gamma)}\ar[d]^{(0,\, 1)}&C\ar@{-->}[r]^{\delta}\ar[d]^{z}&\\
A\ar[r]^{\alpha_A}&X_A\ar[r]^{\beta_A}\ar@{-->}[d]^{0}&A\langle1\rangle\ar@{-->}[r]^{\eta}\ar@{-->}[d]^{\varepsilon}&\\
&&}\end{array}
\end{equation}
of $\E$-triangles.
Since $\alpha_A$ and $\binom{1}{0}$ are $\X$-monic, by Remark \ref{c5}, we have that $y$
is $\X$-monic.
Thus, we have a commutative diagram
\begin{equation}\label{q2}
\begin{array}{l}\xymatrix{B\ar[r]^{y}\ar@{=}[d]&C\ar[r]^z\ar[d]^d&A\langle1\rangle\ar@{-->}[r]^{\varepsilon}\ar[d]^w&\\
B\ar[r]^{\alpha_B}&X_B\ar[r]^{\beta_B}&B\langle1\rangle\ar@{-->}[r]^{\theta}&}\end{array}
\end{equation}
of $\E$-triangles.
It follows that $$\xymatrix{B\ar[r]^{\overline{y}}&C\ar[r]^{\overline{z}\quad}&A\langle1\rangle\ar[r]^{\quad\overline{w}\quad}&B\langle1\rangle}$$
is  a standard triangle in $\cal M$.

It remains to show that $\overline{w}=-\overline{x}\langle1\rangle$. In diagram (\ref{q1}), we obtain a commutative diagram
\begin{equation}\label{q3}
\begin{array}{l}\xymatrix{A\ar[r]^{\alpha_A}\ar[d]^{x}&X_A\ar[r]^{\beta_A}\ar[d]^{-\gamma}&A\langle1\rangle\ar@{-->}[r]^{\eta}\ar[d]^{-1}&\\
B\ar[r]^{y}&C\ar[r]^{z}&A\langle1\rangle\ar@{-->}[r]^{\varepsilon}&}\end{array}
\end{equation}
of $\E$-triangles.
By composing the commutative diagrams (\ref{q2}) and (\ref{q3}), we have the following commutative diagram
$$\xymatrix{A\ar[r]^{\alpha_A}\ar[d]^{x}&X_A\ar[r]^{\beta_A}\ar[d]^{-d\gamma}&A\langle1\rangle\ar@{-->}[r]^{\eta}\ar[d]^{-w}&\\
B\ar[r]^{y}&X_B\ar[r]^{z}&B\langle1\rangle\ar@{-->}[r]^{\theta}&}$$
of $\E$-triangles. It follows that $\overline{w}=-\overline{x}\langle1\rangle$.

(TR3) Suppose that there exists a commutative diagram
$$\xymatrix{A\ar[d]^{\overline{a}}\ar[r]^{\overline{f}}&B\ar[d]^{\overline{b}}\ar[r]^{\overline{g}}&C\ar[r]^{\overline{h}}& A\langle1\rangle\ar[d]^{ \overline{a}\langle1\rangle}\\
A'\ar[r]^{\overline{f'}}&B'\ar[r]^{\overline{g'}}&C'\ar[r]^{\overline{h'}}&
A'\langle1\rangle,}$$
where the rows are standard triangles in $\cal M$.
We can assume that these standard triangles
are induced by
the $\mathbb{E}$-triangles $$\xymatrix{A\ar[r]^{f}&B\ar[r]^{g}&C\ar@{-->}[r]^{\delta}&
}\ \  \textrm{and} \ \ \xymatrix{A'\ar[r]^{f'}&B'\ar[r]^{g'}&C'\ar@{-->}[r]^{\delta'}&
,}$$ where $f,f'$ are $\X$-monic. By the definition of standard triangles of $\cal M$ in Definition \ref{c1}, there exists a (not
necessarily commutative) diagram where the rows are $\mathbb{E}$-triangles in $\C$:
$$\xymatrix{A\ar[d]^{a}\ar[r]^f&B\ar[d]^{b}\ar[r]^g&C\ar@{-->}[r]^{\delta}&\\
A'\ar[r]^{f'}&B'\ar[r]^{g'}&C'\ar@{-->}[r]^{\delta'}&}$$
Since $\overline{b}\circ\overline{f}=\overline{f'}\circ\overline{a}$,  the morphism $bf-f'a$
factors through $\X$, i.e., there exist morphisms $m\colon A\to X$ and
$n\colon X\to B'$ for some $X\in \X$ such that $bf-f'a=nm$. Since
$f$ is $\X$-monic, there exists a morphism $u\colon B\to X$ such that
$m=uf$. Then $bf-f'a=(nu)f$ with $nu\in\X(B,B')$. Set
$b'=b-nu$, then we have $b'f=bf-nuf=f'a$. Then we get the
following commutative diagram of $\mathbb{E}$-triangles in $\C$:
$$\xymatrix{A\ar[d]^{a}\ar[r]^f&B\ar[d]^{b'}\ar[r]^g&C\ar@{-->}[r]^{\delta}&\\
A'\ar[r]^{f'}&B'\ar[r]^{g'}&C'\ar@{-->}[r]^{\delta'}&,}$$
where $f$ and $f'$ are $\X$-monic. Thus the
assertion follows from Lemma \ref{c2}, since $\overline{b'}=\overline{b}$.  \qed

(TR4) Let
$$\xymatrix{A\ar[r]^{\overline{f}}&B\ar[r]^{\overline{f'}}&C\ar[r]&A\langle1\rangle}$$
$$\xymatrix{B\ar[r]^{\overline{g}}&C\ar[r]^{\overline{g'}}&F\ar[r]&B\langle1\rangle}$$
$$\xymatrix{A\ar[r]^{\overline{gf}}&C\ar[r]^{\overline{h'}}&E\ar[r]&A\langle1\rangle}$$
be standard triangles in $\cal M$. We can assume that these standard triangles are induced by the $\E$-triangles
$$\xymatrix{A\ar[r]^{f}&B\ar[r]^{f'}&C\ar@{-->}[r]^{\del_f}&}$$
$$\xymatrix{B\ar[r]^{g}&C\ar[r]^{g'}&F\ar@{-->}[r]^{\del_g}&}$$
$$\xymatrix{A\ar[r]^{gf}&C\ar[r]^{g'}&E\ar@{-->}[r]^{\del_h}&},$$
where $f,g,gf\in\cal A$ are $\X$-monic. By Lemma 3.14 in [NP], we have commutative diagram
$$\xymatrix{
A\ar[r]^{f}\ar@{=}[d]&B\ar[r]^{f'}\ar[d]^{g}&D\ar@{-->}[r]^{\delta_f}\ar[d]^{d}&\\
A\ar[r]^{gf}&C\ar[r]^{h'}\ar[d]^{g'}&E\ar@{-->}[r]^{\delta_h}\ar[d]^{e}&\\
&F\ar@{=}[r]\ar@{-->}[d]^{\del_g}&F\ar@{-->}[d]\\
&&}$$
$$\xymatrix{A\ar[r]^{gf}\ar[d]^f&C\ar[r]^{h'}\ar@{=}[d]&E\ar@{-->}[r]^{\del_h}\ar[d]^e&\\
B\ar[r]^{g}&C\ar[r]^{g'}&F\ar@{-->}[r]^{\del_g}&}$$
of $\E$-triangles.
We will show that $d$ is $\X$-monic. Let $\alpha\colon D\to X$ be a morphism, where $X\in\X$.
Since $g$ is $\X$-monic, there exists a morphism $\beta\colon C\to X$ such that $\beta g=\alpha f'$.
By [NP, Lemma 3.13], there exists a morphism $\gamma\colon E\to X$ which makes the following diagram
commutative.
$$\xymatrix{
  B \ar[d]_{g} \ar[r]^{f'}
                &D  \ar[d]^{d} \ar@/^/[ddr]^{\alpha}  \\
  C \ar[r]_{h'} \ar@/_/[drr]_{\beta}
                &E\ar@{-->}[dr]^{\gamma}\\&& X}$$
It follows that $\alpha=d\gamma$. This shows that $d$ is $\X$-monic.
By Lemma \ref{c2}, we have commutative diagrams
$$\xymatrix{A\ar[r]^{\overline{f}}\ar@{=}[d]&B\ar[r]^{\overline{f'}}\ar[d]^{\overline{g}}&D\ar[r]\ar[d]^{\overline{d}}&A\langle1\rangle\ar@{=}[d]\\
A\ar[r]^{\overline{gf}}&C\ar[r]^{\overline{h'}}\ar[d]^{\overline{g'}}&E\ar[r]\ar[d]^{\overline{e}}&A\langle1\rangle\\
&F\ar@{=}[r]\ar[d]&F\ar[d]\\
&B\langle1\rangle\ar[r]&B\langle1\rangle}$$
$$\xymatrix{A\ar[r]^{\overline{gf}}\ar[d]^{\overline{f}}&C\ar[r]^{\overline{h'}}\ar@{=}[d]&E\ar[r]\ar[d]^{\overline{e}}&A\langle1\rangle\ar[d]^{ \overline{f}\langle1\rangle}\\
B\ar[r]^{\overline{g}}&C\ar[r]^{\overline{g'}}&F\ar[r]&
B\langle1\rangle}$$
of standard triangles in $\cal M$.  \qed

\subsection{Examples}
Since exact categories and extension closed subcategory of a triangulated category are extriangulated categories, we can apply Theorem \ref{h1} to the following situations.

We recall some concepts given by Nakaoka and Palu [NP].
Let $\C$ be an extriangulated category.
\begin{itemize}
\item An object $P\in\C$ is called {\it projective} if
for any $\E$-triangle $\xymatrix{A\ar[r]^{x}&B\ar[r]^{y}&C\ar@{-->}[r]^{\delta}&}$ and any morphism $c\in\C(P,C)$, there exists $b\in\C(P,B)$ satisfying $yb=c$.

We denote the full subcategory of projective objects in $\C$ by $\cal P$. Dually, the full subcategory of injective objects in $\C$ is denoted by $\cal I$.
Remark that $P\in\C$ is projective if and only if $\E(P,C)=0$, for any $C\in\C$.

\item We say $\C$ {\it has enough projectives}, if
for any object $C\in\C$, there exists an $\E$-triangle
$$\xymatrix{A\ar[r]^{x}&P\ar[r]^{y}&C\ar@{-->}[r]^{\delta}&}$$
satisfying $P\in\cal P$.  We can define the notion of having enough injectives dually.

\item $\C$ is said to be {\it Frobenius} if $\C$ has enough projectives and enough injectives
and if moreover the projectives coincide with the injectives.
\end{itemize}

\begin{remark}
In triangulated category $\C$, it is easy to see that $\cal P=\{0\}$ and $\cal I=\{0\}$, $\C$ has enough projectives and enough injectives and $\C$ is a Frobenius extriangulated category.
\end{remark}

\begin{example}\label{lz1}{[NP, Corollary 7.4]}
Let $\C$ be a Frobenius extriangulated category. Then its stable category $\C/\cal I$ is a triangulated category.
\end{example}

\proof It is easy to see that $(\C,\C)$ forms $\cal I$-mutation pair. By Theorem \ref{h1}, we have that
$\C/\cal I$ is a  triangulated category.

\begin{example}\label{lz4}{[IY, Theorem 4.2], [LZ, Theorem 3.11] and [B2, Theorem 3.3]}
Let $\X\subseteq\cal A$ be subcategories of a triangulated category $\C$. If $(\cal A,\cal A)$ forms $\X$-mutation pair and $\cal A$ is extension closed, then the quotient category $\cal A/\X$ is a triangulated category.
\end{example}

\begin{example}\label{lz5}{[B1, Theorem 7.2]}
Let $\C$ be a triangulated category, and $\cal E$ a proper
class of triangles on $\C$. Denote by $\cal I$ the full subcategory of $\C$ consisting of $\cal E$-injective objects, $\cal P$ the full subcategory of $\C$ consisting of $\cal E$-projective objects . If $\C$ has enough $\cal E$-injectives
and enough $\cal E$-projectives and $\cal I=\cal P$, then  the quotient category $\C/\cal I$ is a triangulated category.
\end{example}

\begin{example}{[XZO, Theorem 2.12]}
Let $\X\subseteq\cal A$ be subcategories of an abelian category $\C$. If $\cal A$ is  extension closed and
$(\cal A,\cal A)$ is an $\X$-mutation pair, then the quotient category  $\cal A/\X$  is a triangulated category.
\end{example}

\subsection{One-one correspondence}

\begin{definition}\label{dd1}
Let $\C$ be an  extriangulated category. A subcategory $\X$ of $\C$ is called
\emph{strongly contravariantly finite}, if for any object $C\in\C$, there exists an $\E$-triangle
$$\xymatrix{K\ar[r]&X\ar[r]^{g}&C\ar@{-->}[r]^{\del}&,}$$
where $g$ is a right $\X$-approximation of $C$.

Dually, a subcategory $\X$ of $\C$ is called
\emph{strongly  covariantly  finite}, if for any object $C\in\C$, there exists an $\E$-triangle
$$\xymatrix{C\ar[r]^{f}&X\ar[r]&L\ar@{-->}[r]^{\del'}&,}$$
where $f$ is a left $\X$-approximation of $C$.

A strongly contravariantly finite and strongly  covariantly finite subcategory is called \emph{ strongly functorially finite}.
\end{definition}

\begin{remark}\label{rr1}
Let $\C$ be an extriangulated category, $\X$ a subcategory of $\C$.
\begin{itemize}
\item If $\C$ has enough projectives $\cal P$, then $\X$ is strongly contravariantly finite of $\C$ if and only if $\X$ is contravariantly  finite of $\C$ containing $\cal P$.

\item If $\C$ has enough injectives $\cal I$, then $\X$ is strongly  covariantly finite of $\C$ if and only if $\X$ is  covariantly  finite of $\C$ containing $\cal I$.

\item If $\C$ has enough projectives $\cal P$ and enough injectives $\cal I$, then $\X$ is strongly functorially finite of $\C$ if and only if $\X$ is functorially finite of $\C$ containing $\cal P$ and $\cal I$.

\item If $\C$ is a triangulated category, then $\X$ is strongly contravariantly finite of $\C$ if and only if $\X$ is contravariantly finite of $\C$, $\X$ is strongly covariantly finite of $\C$ if and only if $\X$ is covariantly finite of $\C$, $\X$ is strongly functorially finite of $\C$ if and only if $\X$ is functorially finite of $\C$.

\end{itemize}
\end{remark}
\vspace{3mm}

Let $\X$ be a strongly functorially finite subcategory of an extriangulated category $\C$. Then
for any $A\in\C$, there exist two $\E$-triangles
$$\xymatrix{ A\ar[r]^{x}&X_0\ar[r]&\Sigma A\ar@{-->}[r]^{\delta}&,}$$
and $$\xymatrix{\Omega A\ar[r]&X_1\ar[r]^{f}&A\ar@{-->}[r]^{\eta}&}$$
where $x$ is a left $\X$-approximation and $f$ is a right $\X$-approximation.
We define $\Sigma(A)=\Sigma A$ to be the image of $A$ in $\C/\X$.

For any a morphism $a\colon A\to A'$ in $\C$, there exists a commutative diagram
$$\xymatrix{
A \ar[r]^x \ar[d]_{a} & X\ar[r]^y \ar[d]^{b} & \Sigma A\ar@{-->}[r]^{\delta} \ar[d]^{c}&\\
A'\ar[r]^{x'} & X' \ar[r]^{y'} &\Sigma A' \ar@{-->}[r]^{\delta'}&}$$
of $\E$-triangles in $\C$.
For any $\overline{a}\in\C/\X(A,A')$, define $\Sigma\overline{a}$ to be the image
$\overline{c}$ of $c$ in $\C/\X$.

It is not difficult to see that $\Sigma\colon \C/\X\longrightarrow\C/\X$ is an endofunctor, and
the dual method works for $\Omega$.

For any $\E$-triangle $\xymatrix{ A\ar[r]^{f}&B\ar[r]^g&C\ar@{-->}[r]^{\eta}&}$ in $\C$,
where $f$ is $\X$-monic. Thus we obtain a commutative diagram
$$\xymatrix{
A \ar[r]^f \ar[d]_{u} & B\ar[r]^g \ar[d]^{v} & C\ar@{-->}[r]^{\eta} \ar[d]^{h}&\\
A\ar[r]^{x} & X_0 \ar[r]^{y} &\Sigma A \ar@{-->}[r]^{\delta}&}$$
of $\E$-triangles.

 We define standard right triangles
in $\C/\X$ to be all those isomorphic to the sequences
$$\xymatrix{ A\ar[r]^{\overline{f}}&B\ar[r]^{\overline{g}}&C\ar[r]^{\overline{h}}&\Sigma A}$$
obtained in this way.

Dually, for any $\E$-triangle $\xymatrix@C=0.5cm{ A\ar[r]^{f}&B\ar[r]^g&C\ar@{-->}[r]^{\eta}&,}$ in $\C$,
where $g$ is $\X$-epic.
We define standard left triangles
in $\C/\X$ to be all those isomorphic to the sequences
$\xymatrix@C=0.6cm{ \Omega A\ar[r]^{\overline{f}}&B\ar[r]^{\overline{g}}&C\ar[r]^{\overline{h}}&A}$
obtained in this way.

With this data, the following result holds.
\begin{lemma}{\label{hh1}}{\emph{[BM, Theorem 3.1]}}
Let $\X$ be a strongly functorially finite subcategory of an extriangulated category $\C$.
Then the quotient category $\C/\X$ is a pre-triangulated category.
\end{lemma}

\proof The proof can be carried through in the manner of the proof of [BM, Theorem 3.1].  \qed
\vspace{3mm}

Recall that a full  triangulated subcategory of a pre-triangulated category is called
\emph{thick} if it is closed under direct summands.

\begin{theorem}\label{d4}
If $\X$ is a strongly functorially finite subcategory of an extriangulated category $\C$, then the maps
$$\cal A\longmapsto\cal A/\X \ \ and \ \ \cal N\longmapsto\pi^{-1}(\cal N) $$
give mutually inverse bijection between
\begin{itemize}
\item[\emph{(I)}] $\cal A$ is closed under extensions and $(\cal A,\cal A) $ is $\X$-mutation pair in $\C$.
\item[\emph{(II)}] Thick triangulated subcategories $\cal N$ of $\C/\X$ such that $\pi^{-1}(\cal N)$ is closed under extensions
\end{itemize}
\end{theorem}

\proof If $\cal A$ is closed under extensions and $(\cal A,\cal A) $ is $\X$-mutation pair, by Theorem \ref{h1}, we have that $\cal A/\X$ is triangulated category, and in fact it is a triangulated subcategories of $\C/\X$.
If $\overline{A}=\overline{A_1}\oplus\overline{A_2}$ is a direct sum decomposition in $\C/\X$, where $A\in\cal A$, then
there exist $X_1,X_2\in\X$ such that $A\oplus X_1\simeq A_1\oplus A_2\oplus X_2$ in $\C$.
Since $\cal A$ is closed under direct summands and contains $\X$, we have that $A_1,A_2\in\cal A$ and then
$\overline{A_1},\overline{A_2}\in\cal A/\X$. This shows that $\cal N:=\cal A/\X$ is a thick triangulated subcategories of $\C/\X$.

Conversely, let $\cal N$ be a thick triangulated subcategory of $\C/\X$ such that $\cal A:=\pi^{-1}(\cal N)$ is closed under extensions. Then $\cal A$ is a full subcategory of $\C$ containing $\X$ and clearly $\cal A$ is closed
under direct summands in $\C$. Since $\X$ is a strongly functorially finite subcategory of $\C$ and $\cal N$ is a triangulated subcategory of $\C/\X$, we obtain $\Omega\cal A\subseteq\cal A\supseteq\Sigma\cal A$.

It remains to show that $(\cal A,\cal A)$ is $\X$-mutation pair. For any $A\in\cal A$, since
$\X$ is a strongly functorially finite, there exists an $\E$-triangle
$$\xymatrix{A\ar[r]^f&X_1\ar[r]^{g}&\Sigma A\ar@{-->}[r]^{\del}&,}$$
where $f$ is a left $\X$-approximation of $A$. Since
$\X$ is a strongly functorially finite, there exists an $\E$-triangle
$$\xymatrix{\Omega\Sigma A\ar[r]^u&X_2\ar[r]^{v}&\Sigma A\ar@{-->}[r]^{\del'}&,}$$
where $v$ is a right $\X$-approximation of $\Sigma A$.
Since $\Omega\Sigma A\simeq A$ in $\cal N$, there exists a morphism $a\colon \Omega\Sigma A\to A$ in $\cal A$.
By Definition \ref{a2},
there exists a morphism $b\colon X_2\to X_1$ which makes the following diagram commutative.
$$\xymatrix{
\Omega\Sigma A\ar[r]^u \ar[d]^a & X_2\ar[r]^v \ar@{-->}[d]^{b} & \Sigma A\ar@{-->}[r]^{\delta'}\ar@{=}[d]&\\
A\ar[r]^{f} & X_1 \ar[r]^{g} & \Sigma A\ar@{-->}[r]^{\delta} &}
$$
It follows that $v=gb$.
We claim that $g$ is a right $\X$-approximation of $\Sigma A$. Indeed, for any morphism $\alpha\colon X\to \Sigma A$
where $X\in\X$, since $v$ is a right $\X$-approximation, there exists a morphism $\beta\colon X\to X_2$
such that $v\beta=\alpha$ and then $\alpha=g(bv)$. This shows that $g$ is a right $\X$-approximation.

Similarly, we can show that there exists an $\E$-triangle
$$\xymatrix{\Omega A\ar[r]^x&X_3\ar[r]^{y}&A\ar@{-->}[r]^{\eta}&,}$$
where $x$ is a left $\X$-approximation of $\Omega A$ and $y$ is a right $\X$-approximation of $A$.
Hence $(\cal A,\cal A)$ is $\X$-mutation pair. Clearly the maps
(I) $\longmapsto$ (II) and (II) $\longmapsto$ (I) are mutually inverse.  \qed
\vspace{2mm}

The following corollaries are direct consequences of the above theorem.

\begin{corollary}{\emph{[B2, Theorem 3.3]}} Let $\C$ be a triangulated category.
If $\X$ is a functorially finite subcategory of $\C$, then the maps
$$\cal A\longmapsto\cal A/\X \ \ and \ \ \cal N\longmapsto\pi^{-1}(\cal N) $$
give mutually inverse bijection between
\begin{itemize}
\item[\emph{(I)}] $\cal A$ is closed under extensions and $(\cal A,\cal A) $ is $\X$-mutation pair in $\C$.
\item[\emph{(II)}] Thick triangulated subcategories $\cal N$ of $\C/\X$ such that $\pi^{-1}(\cal N)$ is closed under extensions
\end{itemize}

\end{corollary}

\proof It is easy to see that $\cal A$ is $\X$-Frobenius subcategory if and only if $\cal A$ is extension closed and $(\cal A,\cal A) $ is $\X$-mutation pair. This follows from Theorem \ref{h1} and Theorem \ref{d4}.  \qed

\begin{corollary} Let $\cal B$ be an exact category with enough projectives $\cal P$ and enough injectives $\cal I$.
\begin{itemize}
\item[\emph{(1)}] If $\cal A$ is closed under extensions and $(\cal A,\cal A) $ is $\X$-mutation pair in $\cal B$, then
$\cal A/\X$ is a triangulated category.

\item[\emph{(2)}]If $\X$ is a functorially finite subcategory of $\cal B$ containing $\cal P$ and $\cal I$, then the maps
$$\cal A\longmapsto\cal A/\X \ \ and \ \ \cal N\longmapsto\pi^{-1}(\cal N) $$
give mutually inverse bijection between
\begin{itemize}
\item[\emph{(I)}] $\cal A$ is closed under extensions and $(\cal A,\cal A) $ is $\X$-mutation pair in $\cal B$.
\item[\emph{(II)}] Thick triangulated subcategories $\cal N$ of $\cal B/\X$ such that $\pi^{-1}(\cal N)$ is closed under extensions
\end{itemize}
\end{itemize}
\end{corollary}

\proof This follows from Theorem \ref{h1} and Theorem \ref{d4}.  \qed

\section{J{\o}rgensen theorem}
In this section, we assume that $\C$ is a Krull-Schmidt extriangulated category $k$-category, where $k$ is a field.  Let $f\colon A\to  B$ be a morphism in $\C$.
Recall that from [ASS, ARS] that $f$ is called a \emph{source morphism} if it satisfies
\begin{itemize}
\item $f$ is not a section;

\item If $u\colon A\to C$ is not a section, there exists a morphism $u'\colon B\to C$ such that
$u=u'f$;

\item If $\phi\colon B\to B$ satisfies $\phi f=f$, then $\phi$ is an automorphism.
\end{itemize}
Dually, we can define \emph{a sink morphism}. The notions of a source morphism and a sink morphism
 are also known as \emph{minimal left almost split morphism} and \emph{minimal right split morphism}, respectively.
\vspace{2mm}

Auslander-Reiten theory in a Krull-Schmidt category was established by Liu in [Liu].
However, it can be easily extended to our setting.

\begin{definition}
An $\E$-triangle
$$\xymatrix{A\ar[r]^f&B\ar[r]^{g}&C\ar@{-->}[r]^{\del}&}$$
in $\C$ is called an \emph{Auslander-Reiten $\E$-triangle} if
$f$ is a source morphism and $g$ is a sink morphism.
The translation $\tau$, called the Auslander-Reiten translation, is such that $A=\tau C$ and $C=\tau^{-1}A$ if and only
if $\C$ has an Auslander-Reiten $\E$-triangle  $\xymatrix{A\ar[r]^f&B\ar[r]^{g}&C\ar@{-->}[r]^{\del}&}$.
\end{definition}

\begin{remark}
\  \
\begin{itemize}

\item If $f\colon A\to  B$ is  a source morphism in $\C$, then $A$ is indecomposable. If $g\colon B\to  C$ is a sink morphism in $\C$,  then $C$ is indecomposable.

\item If $\C$ is an abelian category, then an $\E$-triangle
$\xymatrix{A\ar[r]^f&B\ar[r]^{g}&C\ar@{-->}[r]^{\del}&}$
in $\C$ is  an Auslander-Reiten $\E$-triangle if
and only if $\xymatrix@C=0.5cm{0\ar[r]&A\ar[r]^f&B\ar[r]^{g}&C\ar[r]&0}$ is an Auslander-Reiten sequence.
For more details, see [Liu].
\item If $\C$ is a triangulated category, then an $\E$-triangle
$\xymatrix{A\ar[r]^f&B\ar[r]^{g}&C\ar@{-->}[r]^{\del}&}$
in $\C$ is  an Auslander-Reiten $\E$-triangle if
and only if $\xymatrix@C=0.5cm{A\ar[r]^f&B\ar[r]^{g}&C\ar[r]&A[1]}$ is an Auslander-Reiten triangle.
For more details, see [Liu].
\end{itemize}
\end{remark}

We introduce some notations. Let $\X$ be a subcategory of an extriangulated category $\C$. We set
$\underline{\X}:=\X\setminus\cal P$ and $\overline{\X}:=\X\setminus\cal I$. For $A, B\in \C$, $\underline{\Hom}_{\C}(A,B)=\Hom(A,B)/\cal P(A,B)$, $\overline{\Hom}_{\C}(A,B)=\Hom(A,B)/\cal I(A,B)$.

\begin{theorem}\label{jj1}
Let $\C$ be an extriangulated category with Auslander-Reiten translation $\tau$ and $\tau^{-1}$,
$\X$ a strongly functorially finite subcategory of $\C$. Assume that there exists a functorial isomorphism $\E(A,B)\simeq \mathsf{D}\overline{\emph{\Hom}}_{\C}(B,\tau A)\simeq \mathsf{D} \underline{\emph{\Hom}}_{\C}(\tau^{-1}B,A)$, for any $A,B\in\C$, where $\mathsf{D} = \emph{\Hom}_{k}(-,k)$ denotes duality with respect to the base field $k$.
Then the following statements are equivalent:
\begin{itemize}
\item[\emph{(1)}] $(\C,\C)$ is an $\X$-mutation pair.

\item[\emph{(2)}]  The pre-triangulated category $\cal M:=\C/\X$ is a triangulated category.

\item[\emph{(3)}]  $\tau\underline{\X}=\overline{\X}$.

\end{itemize}
\end{theorem}

\proof (1)$\Rightarrow$  (2).   By Theorem \ref{h1}, $\cal M$ is a triangulated category.

(2)$\Rightarrow$ (3).   It suffices to show that $\tau\underline{\X}\subseteq\overline{\X}$
and $\tau^{-1}\overline{\X}\subseteq\underline{\X}$. We only show that $\tau\underline{\X}\subseteq\overline{\X}$, dually, we can show that $\tau^{-1}\overline{\X}\subseteq\underline{\X}$.

For any nonzero indecomposable object $X\in\underline{\X}$, there exists an Auslander-Reiten $\E$-triangle
\begin{equation}\label{t10}
\begin{array}{l}\xymatrix{\tau X\ar[r]^{f}&A\ar[r]^{g}&X\ar@{-->}[r]^{\del}&}\end{array}
\end{equation}
in $\C$.

Since $\tau X$ is not injective object, otherwise, $X=0$. This shows that if $\tau X\in\X$, then $\tau X\in\overline{\X}$, we are done.
 When $\tau X$ is not in $\X$, a morphism in $\C$
from $\tau X$ to an object in $\X$ cannot be a section, so any such morphism factors through $f$, which
is hence an $\X$-monic. This means that there exists a standard triangle
$$
\xymatrix{\tau X\ar[r]^{\overline{f}}&A\ar[r]^{\overline{g}}&X\ar[r]^{h\quad}&\tau X\langle1\rangle}
$$
in $\cal M$, and this is isomorphic to
$$\xymatrix{\tau X\ar[r]^{\overline{f}}&A\ar[r]&0\ar[r]&\tau X\langle1\rangle.}$$
Since $\cal M$ is a triangulated category, we have that $\overline{f}$ is an isomorphism.
Thus there exists an object $X_0\in\X$ such that $A\simeq \tau X\oplus X_0$ in $\C$. Namely
$f$ has the form $\tau X\xrightarrow{\binom{u}{v}}\tau X\oplus X_0$,
and since $X_0\in\X$ and hence zero in $\cal M$, the morphism $\overline{f}=\overline{u}$.
Since $\overline{f}$ is an isomorphism in $\cal M$, then $\overline{u}$ is an isomorphism in $\cal M$,
and so an invertible element of the ring $\textrm{End}_{\cal M}(\tau X)$.
This ring is a quotient of $\textrm{End}_{\C}(\tau X)$, and $\overline{u}$ is the image of $u$.
But $\tau X$ is indecomposable, so $\textrm{End}_{\C}(\tau X)$ is local, and it follows that since $\overline{u}$
is invertible in the quotient, $u$ must itself be invertible.
This implies that $f$ is a section, contradicting that (\ref{t10}) is an Auslander-Reiten $\E$-triangle.
\vspace{2mm}

(3)$\Rightarrow$ (1). Since $\X$ is strongly functorially finite, then for any object $A\in\C$,
there exists an $\E$-triangle in $\C$
$$\xymatrix{A\ar[r]^{f}&X_0\ar[r]^{g}&B\ar@{-->}[r]^{\del}&,}$$
where $f$ is a left $\X$-approximation of $A$. We will prove that $g$ is a right $\X$-approximation. Applying the functor $\Hom_{\C}(\X,-)$ to the
above  $\E$-triangle, we obtain an exact sequence
$$\Hom_{\C}(\X,X_0)\xrightarrow{\Hom_{\C}(\X,g)}\Hom_{\C}(\X,B)\xrightarrow{\ \ \ }
\E(\X,A)\xrightarrow{\E(\X,f)}\E(\X,X_0).$$
Thus we know that  $$\Hom_{\C}(\X,X_0)\xrightarrow{\Hom_{\C}(\X,g)}\Hom_{\C}(\X,B)$$
is an epimorphism if and only if
$$\E(\X,A)\xrightarrow{\E(\X,f)}\E(\X,X_0)$$
is a monomorphism if and only if
$$\E(\underline{\X},A)\xrightarrow{\E(\underline{\X},f)}\E(\underline{\X},X_0)$$
is a monomorphism.
By functorial isomorphism, it suffices to show that the morphism
$$\overline{\Hom}_{\C}(X_0,\tau\underline{\X})\xrightarrow{\Hom_{\C}(f,\;\tau\underline{\X})}\overline{\Hom}_{\C}(A,\tau\underline{\X})$$
is an epimorphism.
That is to say, $$\overline{\Hom}_{\C}(X_0,\overline{\X})\xrightarrow{\Hom_{\C}(f,\; \overline{\X})}\overline{\Hom}_{\C}(A,\overline{\X})$$
is an epimorphism.
For this, observe that we have a commutative diagram
$$\xymatrix@C=2.2cm{\Hom_{\C}(X_0,\X)\ar[r]^{\Hom_{\C}(f,\, \X)}\ar[d]&\Hom_{\C}(A,\X)\ar[d]\\
\overline{\Hom}_{\C}(X_0,\overline{\X})\ar[r]^{\overline{\Hom}_{\C}(f,\;\overline{\X})}&\overline{\Hom}_{\C}(A,\overline{\X}),}$$
where the vertical morphisms are natural epimorphisms. Hence it is enough to show that the
morphism
$$\Hom_{\C}(X_0,\X)\xrightarrow{\Hom_{\C}(f,\; \X)}\Hom_{\C}(A,\X)$$
is an epimorphism. This shows that $g$ is a right $\X$-approximation of $B$.

Dually, we can show that for any $B'\in\C$, there exists an exact sequence
$$\xymatrix{A\ar[r]^{f'}&X'\ar[r]^{g'}&B'\ar@{-->}[r]^{\del'}&,}$$
where $g'$ is a right $\X$-approximation and $f'$ is a left $\X$-approximation.

Therefore, $(\C,\C)$ is an $\X$-mutation pair.  \qed
\vspace{2mm}

As applications of Theorem \ref{jj1} to triangulated categories or to abelian categories, we obtain J{\o}rgensen's Theorem 3.3 in [J] and its abelian version.

\begin{corollary}\label{cor1}{\emph{[J, Theorem 3.3]}}
Let $\C$ be a triangulated category with Auslander-Reiten translation $\tau$, $\X$ a functorially finite subcategory of $\C$. Then
the pre-triangulated category $\C/\X$ is a triangulated category  if and only if $\tau\X=\X$.
\end{corollary}

\proof This follows from Theorem \ref{jj1} and Remark \ref{rr1}. \qed

\begin{corollary}\label{cor2}
Let $\C=\emph{mod}\Lambda$, where $\Lambda$ is an Artin algebra.
Assume that $\X$ is a functorially finite subcategory of $\C$ containing $\cal P$ and $\cal I$. Then the pre-triangulated category
$\C/\X$ is a triangulated category  if and only if $\tau\underline{\X}=\overline{\X}$.
\end{corollary}

\proof This follows from Theorem \ref{jj1} and Remark \ref{rr1}.  \qed
\vspace{4mm}

There are many examples suit for conditions in Corollary \ref{cor1},  see [J], we give an example for Corollary \ref{cor2},

\begin{example}
Consider the finite dimensional $k$-algebra, where $k$ is a field.
$$\Lambda:=kQ/R,$$
 where $Q$ is the following quiver
$$\xymatrix@C=1.5cm{ \ar@(ul,dl)_{\alpha}1\dr^{\beta} & 2 \ar@(ur,dr)^{\gamma}}$$
and $R$ is generated by $\alpha^2,\gamma^2$ and $\beta\alpha-\gamma\beta$.

Take $\X:=\add(P_1\oplus P_2\oplus I_1)\subseteq \textrm{mod}\Lambda:=\C$.
We claim that $\tau \underline{\X}=\overline{\X}$, where $\tau$ is the Auslander-Reiten translation.
Since $P_1 =I_2$ is a projective and injective object,
it suffices to show that
$\tau I_1=P_2$. Indeed, take the minimal projective resolution of $I_1$
$$\xymatrix{0\ar[r]&P_2\ar[r]&P_1\ar[r]&I_1\ar[r]&0.}$$
It is easy to see that $\tau I_1=P_2$.
By Corollary \ref{cor2}, we have that $\C/\X$ is a triangulated category.
\end{example}

\begin{corollary}\label{cc1}
Under the same assumptions as in Theorem \emph{\ref{jj1}}. Let
$$\xymatrix{A\ar[r]^f&B\ar[r]^{g}&C\ar@{-->}[r]^{\del}&}$$
be an $\E$-triangle in $\C$. If $\tau\underline{\X}=\overline{\X}$, then
$f$ is $\X$-monic if and only if $g$ is $\X$-epic.
\end{corollary}

\begin{theorem}\label{ppp1}
Let $(\C,\E,\s)$ be an extriangulated category with Auslander-Reiten translation $\tau$ and $\tau^{-1}$,
$\X$ a strongly functorially finite subcategory of $\C$ and satisfying $\tau\underline{\X}=\overline{\X}$. Assume that there exists a functorial isomorphism $\E(A,B)\simeq \mathsf{D}\overline{\emph{\Hom}}_{\C}(B,\tau A)\simeq \mathsf{D} \underline{\emph{\Hom}}_{\C}(\tau^{-1}B,A)$, for any $A,B\in\C$.
For any $A,C\in\C$, define $\E'(C,A)$ to be the
collection of $\E$-triangles of the form
$\xymatrix{A\ar[r]^f&B\ar[r]^{g}&C\ar@{-->}[r]^{\del}&}$, where $f$ is $\X$-monic.
$\s'(\delta)=[A\xrightarrow{~f~}B\xrightarrow{~g~}C]$, for any $\del\in\E'(C,A)$.  Then $(\C,\E',\s')$ is a Frobenius extriangulated category.
\end{theorem}

\proof Firstly, we show that $(\C,\E',\s')$ is an extriangulated category.

(ET1) Let $\xymatrix@C=0.7cm{A\ar[r]^f&B\ar[r]^{g}&C\ar@{-->}[r]^{\del}&}$
and $\xymatrix@C=0.7cm{A'\ar[r]^{f'}&B'\ar[r]^{g'}&C'\ar@{-->}[r]^{\del'}&}$ be $\E$-triangles in $(\C,\E,\s)$, where $f,f'$ are $\X$-monic. Since $\s$ is an additive realization of $\E$, there exists a morphism $y\colon B\to B'$ which gives the following morphism of $\E$-triangles
$$\xymatrix{A\ar[r]^{f}\ar@{=}[d]&B\ar[r]^g\ar@{-->}[d]^y&C\ar@{-->}[r]^{\delta}\ar@{=}[d]&\\
A\ar[r]^{f'}&B'\ar[r]^{g'}&C\ar@{-->}[r]^{\delta'}&.}$$
By Lemma \ref{c6}, we have that $y$ is an isomorphism. This shows that $s'(\del)$ is the equivalence class.
Thus $\E'(C,A)$ is the
collection of all equivalence classes of $\E$-triangles of the form
$\xymatrix{A\ar[r]^f&B\ar[r]^{g}&C\ar@{-->}[r]^{\del}&}$, where $f$ is $\X$-monic.
We obtain a biadditive functor $\E'\colon \C\op\times\C\to Ab$. Its structure is given as follows.
\begin{itemize}
\item For any $\delta\in\E'(C,A)$ and any morphism $a\in\C(A,A')$, we can get the following morphism
of $\E$-triangles
$$\xymatrix{A\ar[r]^{f}\ar[d]^{a}&B\ar[r]^g\ar[d]^b&C\ar@{-->}[r]^{\delta}\ar@{=}[d]&\\
A'\ar[r]^{f'}&B'\ar[r]^{g'}&C\ar@{-->}[r]^{\delta'}&.}$$
By Corollary \ref{cc1}, we have that $g$ is $\X$-epic. We will show that $g'$ is $\X$-epic.
For any morphism $\alpha\colon X\to C$, where $X\in\X$, since $g$ is $\X$-epic, there exists a morphism
$\beta\colon X\to B$ such that $\alpha=g\beta$. It follows that $\alpha=g'(b\beta)$.
This shows that $g'$ is $\X$-epic. By Corollary \ref{cc1}, we have that $f'$ is $\X$-monic.
This gives $\E'(C,a)(\del)=a_\ast\del=[A'\xrightarrow{~f'~}B'\xrightarrow{~g'~}C]$, where $f'$ is $\X$-monic.
Dually, for any $c\in\C(C',C)$, the morphism $\E'(c,A)=c^{\ast}\colon
\E'(C,A)\to \E'(C',A)$ is defined.
\item The zero element in $\E'(C,A)$ is given by the split $\E$-triangle
$0=[A\xrightarrow{\binom{1}{0}}A\oplus C\xrightarrow{(0,\  1)}C].$
For any $\delta,\del'\in\E'(C,A)$, its sum $\del+\del'$ is given by
$\delta+\delta'=\E'(\Delta_C,\nabla_A)[\delta\oplus\del'],$
where $\Delta_C=\binom{1}{1}\colon C\to C\oplus C,\  \nabla_A=(1,1)\colon A\oplus A\to  A$.
It is straightforward to verify that $\E'(C,A)$ is additive group.
This shows that $\E'$ is biadditive functor.
\end{itemize}

(ET2) $s'$ is  a correspondence which associates an equivalence class
$\s'(\delta)=[A\xrightarrow{~f~}B\xrightarrow{~g~}C]$
to any $\E'$-extension $\delta\in\E'(C,A)$. Then (ET2) is trivially satisfied.
\vspace{2mm}

(ET3) (ET3$\op$) are trivially satisfied.

(ET4) Let $\xymatrix@C=0.7cm{A\ar[r]^{f} & B \ar[r]^{f'} & D\ar@{-->}[r]^{\delta}&}$ and $ \xymatrix@C=0.7cm{B\ar[r]^{g} & C \ar[r]^{g'} & F\ar@{-->}[r]^{\delta'}&}$
be $\E$-triangles  in $(\C,\E,\s)$, where $f,g$ are $\X$-monic. Then by (ET4) for $(\C,\E,\s)$
there exist an object $E\in\C$, a commutative diagram
$$\xymatrix{
A\ar[r]^{f}\ar@{=}[d]&B\ar[r]^{f'}\ar[d]^{g}&D\ar@{-->}[r]^{\delta}\ar[d]^{d}&\\
A\ar[r]^{gf}&C\ar[r]^{h'}\ar[d]^{g'}&E\ar@{-->}[r]^{\delta''}\ar[d]^{e}&\\
&F\ar@{=}[r]\ar@{-->}[d]^{\del}&F\ar@{-->}[d]\\
&&}$$
in $(\C,\E,\s)$, and an $\E$-extension $\del^{''}\in\E(E,A)$ realized by $\xymatrix@C=0.7cm{A\ar[r]^{h} & C \ar[r]^{h'} & E},$ which satisfy the following compatibilities.
\begin{itemize}
\item[{\rm (i)}] $\xymatrix@C=0.7cm{D\ar[r]^{d} & E \ar[r]^{e} & F}$  realizes $f'_{\ast}\del'$,
\item[{\rm (ii)}] $d^\ast\del''=\del$,
\item[{\rm (iii)}] $f_{\ast}\del''=e^{\ast}\del'$.
\end{itemize}
It suffices to show that $gf$ is $\X$-monic and $d$ is $\X$-monic.
Since $f,g$ are $\X$-monic, then $gf$ is $\X$-monic.
Now we will show that $d$ is $\X$-monic. By Corollary \ref{cc1}, we only show that
$e$ is $\X$-epic. Let $\alpha\colon X\to F$, where $X\in\X$. By Corollary \ref{cc1}, we have that
$g'$ is $\X$-epic. So there exists a morphism $\beta\colon X\to C$ such that $\alpha=g'\beta$.
Since $g'=eh'$, we have $\alpha=e(h'\beta)$. This shows that $d$ is $\X$-monic.
(ET4$\op$) can be shown dually.
\vspace{2mm}

Secondly, we show that $(\C,\E',\s')$ is Frobenius.
By the definition of $s'$ and Corollary \ref{cc1}, every object in $\X$ is projective and injective in $(\C,\E',\s')$.
Since $\X$ is a strongly functorially finite subcategory of $\C$, we have that
$(\C,\E',\s')$ has enough projectives and enough injectives.
Take any projective object $P$ in $(\C,\E',\s')$, since $\X$ is a strongly functorially finite subcategory of $\C$, there exists an $\E$-triangle
$$\xymatrix{K\ar[r]^x&X\ar[r]^{y}&P\ar@{-->}[r]^{\eta}&,}$$
where $y$ is a right $\X$-approximation of $P$.
 Thus this $\E$-triangle splits implies that $P\in\X$. The injective objects can be treated similarly.
 Therefore, the classes of projective objects and injective objects in $(\C,\E',\s')$
 both equal $\X$.
 This completes the proof. \qed

\begin{corollary}{\emph{[KIWY, Proposition 2.16]}}
Let $(\cal B,\cal S)$ be an exact category with Auslander-Reiten translation $\tau$ and $\tau^{-1}$
and with enough projectives $\cal P$ and enough injectives $\cal I$,
$\X$ a functorially finite subcategory of $\cal B$ containing $\cal P$ and $\cal I$, which satisfies $\tau\underline{\X}=\overline{\X}$. Assume that there exists a functorial isomorphism $\E(A,B)\simeq \mathsf{D}\overline{\emph{\Hom}}_{\C}(B,\tau A)\simeq \mathsf{D} \underline{\emph{\Hom}}_{\C}(\tau^{-1}B,A)$, for any $A,B\in\C$.
Then the following statements hold.
\begin{itemize}
\item[\emph{(1)}]  Let $\xymatrix@C=0.5cm{0\ar[r]&A\ar[r]^f&B\ar[r]^{g}&C\ar[r]&0}$ be an exact sequence in $(\cal B,\cal S)$.
Then $f$ is $\X$-monic if and only if $g$ is $\X$-epic. We denote the subclass of the kernel-cokernel pair in
$\cal S$ satisfying one of these equivalent conditions by $\cal S'$.

\item[\emph{(2)}] $(\cal B,\cal S')$ is a Frobenius category whose projective objects are precisely $\X$.
\end{itemize}
\end{corollary}

\proof An exact category  can
be regarded as an extriangulated category, whose inflations are monomorphic and
whose deflations are epimorphic. Conversely, if an extriangulated category $\C$, in which any inflation
is monomorphic, and any deflation is epimorphic, then $\C$ is an exact category. See [NP, Corollary 3.18].
This follows from Corollary \ref{cc1} and Proposition \ref{ppp1}.  \qed

\begin{corollary}\label{cc2}
Let $\C$ be a triangulated category with Auslander-Reiten translation $\tau$, and
$\X$ a functorially finite subcategory of $\C$, which satisfies $\tau\X=\X$.
Then the following statements hold.
\begin{itemize}
\item[\emph{(1)}]  Let $\xymatrix@C=0.5cm{A\ar[r]^f&B\ar[r]^{g}&C\ar[r]&A[1]}$ be an triangle in $\C$.
Then $f$ is $\X$-monic if and only if $g$ is $\X$-epic.

\item[\emph{(2)}] For any $A,C\in\C$, define $\E'(C,A):=\C(C,A[1])$ to be the
collection of all equivalence classes of triangles of the form
$\xymatrix@C=0.7cm{A\ar[r]^f&B\ar[r]^{g}&C\ar[r]^{\del\;}&A[1]}$, where $f$ is $\X$-monic.
$\s'(\delta)=[A\xrightarrow{~f~}B\xrightarrow{~g~}C]$, for any $\del\in\E'(C,A)$.  Then $(\C,\E',\s')$ is a Frobenius extriangulated category whose projective objects are precisely $\X$.
\end{itemize}
\end{corollary}

\proof This follows from Corollary \ref{cc1} and Proposition \ref{ppp1}.  \qed

\begin{remark}\label{rem1} When $\X\not= \{0\}$, the new Frobenius extriangulated category $(\C,\E',\s')$ in Corollary above is not triangulated, since $\X$ is projectives and injectives and non-zero. When $\X\not= \C$, it is easy to see that the new Frobenius extriangulated category $(\C,\E',\s')$ is not exact. Otherwise any $\E'-$extension splits and then any object in $\C$ is projective and injective. Then $\X=\C$, a contradiction. We take a subcategory $\{0\}\not=\X\varsubsetneqq \C$ satisfying the condition in Corollary above, Then the new Frobenius extriangulated category $(\C,\E',\s')$ is neither exact nor triangulated. This provides many new extriangulated categories which are not exact or triangulated, see [NP].

\end{remark}

Panyue Zhou\\
Department of Mathematical Sciences, Tsinghua University,
100084 Beijing, P. R. China.\\
E-mail: \verb"pyzhou@math.tsinghua.edu.cn"\\[0.3cm]
Bin Zhu\\
Department of Mathematical Sciences, Tsinghua University,
100084 Beijing, P. R. China.\\
E-mail: \verb"bzhu@math.tsinghua.edu.cn"


\begin{thebibliography}{99}
\bibitem[AR]{ar} M. Auslander, I. Reiten. Applications of contravariantly finite subcategories. Adv. Math., 86(1), 111-152, 1991.




\bibitem[ARS]{ars} M. Auslander, I. Reiten, S. O. Smal{\o}. Representation theory of Artin algebras. Cambridge Studies in Advanced Mathematics, 36. Cambridge University Press, Cambridge, 1997.


\bibitem[ASS]{ass} I. Assem, D. Simson, A. Skowronski. Elements of the Representation Theory of Associative Algebras. Volume 1: Techniques of Representation Theory. Cambridge University Press, 2006.




\bibitem[B1]{b1} A. Beligiannis. Relative homological algebra and purity in triangulated categories. J. Algebra, 227(1): 268-361, 2000.


\bibitem[B2]{b2} A. Beligiannis. Rigid objects, triangulated subfactors and abelian localizations. Math. Z., 274(3), 841-883, 2013.


\bibitem[BMRRT]{bmrrt} A. Buan, R. Marsh, M. Reineke, I. Reiten, G. Todorov. Tilting theory and cluster combinatorics. Adv. Math. 204(2), 572-618, 2006.



\bibitem[B\"{u}]{bu} T. B\"{u}hler. Exact categories. Expo. Math., 28(1), 1-69, 2010.



\bibitem[BM]{bm}  A. Beligiannis, N. Marmaridis. Left triangulated categories arising from contravariantly finite subcategories. Comm. Algebra,  22(12), 5021-5036, 1994.



\bibitem[FZ]{fz}
S. Fomin, A. Zelevinsky. Cluster algebras. I. Foundations.
 J. Amer. Math. Soc. 15(2), 497-529, 2002.


\bibitem[GLS]{gls}
C. Gei{\ss}, B. Leclerc, J. Shr\"{o}er. Rigid modules over preprojective algebras.
 Invent. math. 165(3), 589-632, 2006.

\bibitem[Ha]{ha} D. Happel. Triangulated categories in the representation of finite dimensional algebras.  Vol. 119, Cambridge University Press, 1988.


\bibitem[IY]{iy} O. Iyama, Y. Yoshino. Mutations in triangulated categories and rigid Cohen-Macaulay modules. Invent. Math., 172(1), 117-168, 2008.

\bibitem[J]{j} P. J{\o}rgensen. Quotients of cluster categories. Proc. Roy. Soc. Edinburgh Sect. A,  140(1), 65-81, 2010.

\bibitem[Ke]{ke}B. Keller.
Derived categories and their uses. Handbook of algebra, Vol. 1, 671-701, North-Holland, Amsterdam, 1996.


\bibitem[KIWY]{ke} M. Kalck, O. Iyama, M. Wemyss, D. Yang. Frobenius categories, Gorenstein algebras and rational surface singularities, Compos. Math., 151(3), 502-534, 2015.






\bibitem[Li]{li} Z.W. Li. A homotopy theory of additive categories with suspensions. arXiv:1510.02258v2, 2016.

\bibitem[Liu]{liu} S. P. Liu.
Auslander-Reiten theory in a Krull-Schmidt category.
S\~{a}o Paulo J. Math. Sci., 4(3), 425-472, 2010.



\bibitem[LZ]{lz} Y. Liu, B. Zhu. Triangulated quotient categories. Comm. Algebra, 41(10), 3720-3738, 2013.




\bibitem[NP]{np}H. Nakaoka, Y. Palu. Mutation via Hovey twin cotorsion pairs and model structures in extriangulated categories.  arXiv: 1605.05607, 2016.









\bibitem[XZO]{xzo} J. D. Xu, P. Y. Zhou, B. Y. Ouyang. Mutation pairs in abelian categories.
Comm. Algebra, 44(7),  2732-2746, 2016.







\end{thebibliography}
\end{document}